\documentclass[10pt]{article}
\usepackage{amssymb}
\usepackage{amsmath}
\usepackage{amsfonts}
\usepackage{amscd}
\usepackage{amsthm}
\usepackage{graphicx}
\usepackage{enumerate}
\usepackage{empheq}
\usepackage{xcolor}
\usepackage{subcaption}
\usepackage{cancel}
\usepackage{multirow}
\usepackage{cite}
\usepackage{float}
\usepackage[colorlinks=true]{hyperref}
\hypersetup{urlcolor=blue, citecolor=blue, linkcolor=blue}

\numberwithin{equation}{section}
\newtheorem{theorem}{Theorem}[section]

\newtheorem{lemma}{Lemma}[section]
\newtheorem{example}{Example}[section]

\newcommand{\bn}{\textbf{n}}
\newcommand{\bx}{\textbf{x}}

\DeclareMathOperator*{\esssup}{ess\: sup}

\theoremstyle{remark}
\newtheorem{remark}{Remark}[section]
\newcommand{\vertiii}[1]{{\left\vert\kern-0.25ex\left\vert\kern-0.25ex\left\vert #1 
		\right\vert\kern-0.25ex\right\vert\kern-0.25ex\right\vert}}
\begin{document}
\title{
Discontinuous Galerkin Approximation to a Class of  Strongly Damped Wave Equations with Uniform Exponential Decay}
\author{%
\begin{tabular}{@{}c@{}}
Debendra K. Swain\textsuperscript{1},
P. Danumjaya\textsuperscript{2},
Anil Kumar\textsuperscript{2}
and Amiya K. Pani\textsuperscript{2}\\[4pt]
\small\textsuperscript{1}Tata Institute of Fundamental Research, Center for Applicable Mathematics,\\
\small Bengaluru--560065, India.
\\
\small\textsuperscript{2}Department of Mathematics,
Birla Institute of Technology and Science, Pilani\\
\small K K Birla Goa Campus, Zuarinagar, Sancoale,
Goa--403726, India.\\[4pt]
\small Email: 1729etal@gmail.com, danu@goa.bits-pilani.ac.in;\\
\small anilpundir@goa.bits-pilani.ac.in
and amiyap@goa.bits-pilani.ac.in
\end{tabular}%
}
\date{}
	\maketitle
\begin{abstract}
This paper proposes and examines a completely discrete discontinuous Galerkin finite element scheme to approximate the solution of  a class of strongly damped wave equations with uniform exponential decay property.
In the first part, a semidiscrete discontinuous Galerkin approximation is employed in the spatial variable, while keeping time variable  continuous
and optimal error estimates are established while preserving the exponential decay property of the original system. Subsequently, a fully discrete scheme is developed and analyzed. Using suitable energy techniques, the consistency between the decay properties of the continuous and discrete systems is demonstrated. In particular, the constants appearing in various estimates remain uniformly bounded even as the damping parameters approach zero. Moreover, error estimates that preserve the uniform decay property are highlights of this paper.
Further extension of our analysis includes a class of strongly damped fourth  order evolution equations using SIPG as well as $C^0$-IP methods.
Finally, several numerical experiments are performed to confirm the theoretical convergence results, verify the predicted decay behavior, and investigate the influence of model parameters on stability and computational accuracy.
\end{abstract}
\vspace{1em} 
	
	\section{Introduction}\label{s1}
This work is devoted to the discontinuous Galerkin error analysis and asymptotic behavior of solutions to a class of damped wave equations of the following type. Given a bounded domain $\Omega \subset \mathbb{R}^2$  with boundary $\partial\Omega,$ seek
$u=u(x,t)$ for $(x,t)\in \Omega \times [0,\infty)$ such that
\begin{eqnarray}\label{1.1}
u'' - \beta \Delta u' + \alpha u' -\Delta u = 0 \quad \text{in } \Omega \times (0,\infty) 
\end{eqnarray}
 with boundary condition
\begin{eqnarray}
u(x,t) = 0 \quad \text{on } \partial\Omega \times [0,\infty), \label{1.2}
\end{eqnarray}
and initial  conditions
\begin{eqnarray}
u(x,0)=u_0(x), \qquad u'(x,0)=u_1(x) \quad \text{in } \Omega. \label{1.3}
\end{eqnarray}
Here, $\Delta$ denotes the Laplacian, and
$u':=\frac{\partial u}{\partial t}$, $u''=\frac{\partial^2 u}{\partial t^2}$. The parameters $\alpha$ and $\beta$ represent damping coefficients and satisfy one of the following conditions: $(i)$ $\alpha>0$ and $\beta>0$, $(ii)$
 $\alpha>0$ and $\beta\ge 0$,
$(iii)$ $\alpha\ge 0$ and $\beta>0$.
When $\alpha>0$ and $\beta=0$, equation \eqref{1.1} corresponds to the weakly damped wave equation. In contrast, the case $\alpha=0$ and $\beta>0$ gives rise to the strongly damped wave equation. Such equations, together with their nonlinear counterparts, appear in a wide range of scientific and engineering applications, including vibration suppression, acoustic wave propagation, viscoelasticity, railway track dynamics, and quantum-mechanical models \cite{DellOro2011, Edalatzadeh2018, Ghidaglia1991, Graber2016}. Their mathematical analysis, particularly concerning the well-posedness, stability, and long-time asymptotic behavior has attracted considerable attention in the literature, see, for example, \cite{chen1990, chen1979, chen1981,Russell1975, rauch1976,Temam2012, Massatt1983, AngDinh1988, Carvalho2008, Danumjaya2025, Pata2005, PataZelik2006} and references therein. To be more specific, the long-time behavior of weakly damped wave equations (\(\beta=0\)) has attracted considerable attention. Explicit non-uniform decay estimates have been established through control-theoretic methods and have subsequently been refined using the spectral properties of the operator \(A\) and variable damping coefficients. Further developments and related results have been reported in \cite{Russell1975, Temam2012, rauch1976, chen1979, chen1990, chen1981}.

For conforming finite element approximations of \eqref{1.1}--\eqref{1.3}, we refer to \cite{Danumjaya2025, larsson1991, rincon2013, egger2019, ljung2021, zhao2023} and  references therein. Especially, Larsson {\it et al.} \cite{larsson1991} have established optimal error estimates for finite element discretizations under minimal regularity assumptions on the initial data, with convergence results remaining uniform in time. For weakly damped wave equations, T\'ebou {\it et al.} \cite{tebou2003} proposed a semidiscrete finite difference approximation and investigated the uniform  decay properties of the associated discrete energy by introducing an artificial numerical viscosity term, see also \cite{munch2007}. Rincon {\it et al.} \cite{rincon2013}  have analyzed finite element approximations of a locally damped wave model, proving well-posedness and deriving error estimates which are not uniform in time. However, the asymptotic stability of the discrete scheme has not been  addressed in their article. Later, Egger {\it et al.} \cite{egger2019} have employed mixed finite element technique to a one dimensional weakly internal damped wave  equation,which  preserves uniform exponential decay at both primary and higher energy levels. More recently, Danumjaya {\it et al.} \cite{Danumjaya2025} have demonstrated that conforming finite element approximations inherit the exponential stability of weakly damped wave equations while achieving optimal convergence rates. Kumar {\it et al.} \cite{Krishan2025} have investigated the asymptotic behavior and finite element error analysis under various damping regimes and established exponential decay estimates for the continuous, semidiscrete, and fully discrete solutions, showing that the numerical approximations preserve the decay properties. The resulting error estimates remain uniform as the damping parameters tend to zero.

The numerical analysis available in  literature for \eqref{1.1}--\eqref{1.3} has primarily been carried out using conforming finite element techniques, while there are a few articles related to the development and analysis of discontinuous Galerkin (DG) methods. The DG framework offers several attractive features for the numerical approximation of partial differential equations. In particular, it naturally accommodates discontinuities across element interfaces, leads to block-diagonal mass matrices that facilitate efficient computations, and is well-suited for parallel implementations. Moreover, DG methods can easily handle nonconforming meshes and variable polynomial degrees, providing considerable flexibility in practical computations. Their local conservation properties, robustness, and high-order accuracy make them especially effective for problems involving complex solution structures. These advantages motivate the investigation of DG approximations for the damped wave problem considered in this work which preserves the uniform exponential decay property. Discontinuous Galerkin methods have been investigated for
both spatial and temporal approximations of wave problems.
Grote and Mitkova~\cite{GroteMitkova2013} have developed
high-order explicit local time-stepping schemes for damped
wave equations based on Adams-Bashforth methods, allowing
continuous and discontinuous Galerkin spatial approximations
on locally refined meshes. Subsequently, Antonietti
{\it et al.}~\cite{AntoniettiEtAl2025Review} have reviewed
DG time discretizations for wave propagation problems,
comparing first- and second-order formulations with
respect to stability, accuracy, and computational efficiency.

This paper aims to develop and analyze a DG formulation that preserves the intrinsic dissipative structure of the underlying model. We derive explicit exponential decay behavior  for the error estimates of both semidiscrete and fully discrete approximations. The derived decay rates are expressed in terms of the damping coefficients and the coercive constant of the associated discrete elliptic operator. Our analysis further shows that both  semidiscrete and fully discrete DG schemes inherit the stability and long-time decay properties of the continuous system.

The main contributions of this work are  summarized as follows.
\begin{itemize}
\item We develop a semidiscrete discontinuous Galerkin approximation for the strongly damped wave equation and prove that the resulting scheme inherits the exponential stability of the continuous problem. 
Further, optimal-order error estimates preserving exponential decay behavior  are established in both \(L^{\infty}(L^{2}(\Omega))\) and energy--norms.

\item By combining the spatial DG discretization with a finite-difference approximation in time, we develop a fully discrete scheme and  derive exponential decay estimates that remain uniform with respect to the spatial and temporal discretizing parameters.
\item The analysis is further extended to a fourth order damped wave equation using both SIPG and $C^0$-interior penalty methods, which preserve uniform exponential decay properties. 

\item To validate our theoretical results, a series of numerical experiments is presented. The computations confirm the predicted convergence rates, demonstrate the preservation of long-time decay behavior, and illustrate the influence of  damping parameters on the stability and accuracy of computed solutions.
\end{itemize}

The paper is organised as follows. Section~\ref{s2} introduces
the mathematical setting and weak formulation, together with
the continuous decay and regularity estimates used in the analysis.
Section~\ref{s3} develops the semidiscrete DG approximation
and establishes its exponential decay properties and error estimates.
Section~\ref{s4} presents the fully discrete scheme and derives
the corresponding decay and error bounds.
Section~\ref{s5} extends the analysis to nonhomogeneous problems
and fourth-order evolution equations.
Section~\ref{s6} presents numerical experiments that assess
the predicted convergence rates and illustrate the long-time
behaviour of the numerical solutions.
Finally, Section~\ref{s7} summarises the main findings.
\section{Preliminaries and Weak Formulation}\label{s2}
We begin this section by describing the spaces. Now, consider the standard Sobolev spaces $W^{r,p}(\Omega)$ (\cite{adams9247sobolev}) on a domain $\Omega\subset\mathbb{R}^2$:
\[W^{r,p}(\Omega)=\{v\in L^p(\Omega):\, \partial^\beta v\in L^p(\Omega), \; \forall \; |\beta|\leq r\},\]
where $r\geq 0$ is an integer,  $ 1\leq p \leq \infty$
and  $\beta =(\beta_1,\beta_2)$ is a multi-index with non-negative integers $\beta_i$, $i=1,2$, where \(\displaystyle{|\beta|= \sum_{i=1}^{2} \beta_i}\) and $\partial^\beta v$ represents the partial derivative of order $r.$ 
For $H^r(\Omega):=W^{r,2}(\Omega)$, the norm and semi-norm are denoted by $\|\cdot\|_{r,\Omega}$ and $|\cdot|_{r,\Omega}$, respectively. Whenever there is no confusion, we  represent the norm and semi-norm in $H^r(\Omega)$ as $\|\cdot\|_{r}$ and $|\cdot|_{r}$, respectively.
The notation $\|\cdot\|_{L^p(\Omega)}$ denotes the norm in the Lebesgue space $L^p(\Omega)$.  For $p = 2$, we simply write  the norm and inner-product in $L^2(\Omega)$ as  $\|\cdot\|$ and $(\cdot,\cdot)$, respectively.
Further, $H_0^1(\Omega)$ is the space of functions in $H^1(\Omega)$ with zero traces on the boundary $\partial \Omega$.
The notation $L^p(0, T; V)$ represents the space of measurable functions $\phi$, such that,  for each $t \in(0,T)$, $\phi(t)$ takes values in a Banach space $V$, satisfying
\[ \int_0^T \|\phi(t)\|^p_V~dt <\infty~~\mbox {for}~1 \le p < \infty, ~~~\mbox{and}~~~
{\esssup_{0<t<T}}\: \|\phi(t)\|_V <\infty~~\mbox {for}~p=\infty. \]

The variational formulation of  the problem \eqref{1.1}--\eqref{1.2} is to find   
\(
u:(0,\infty)\to H_0^1(\Omega)
\)
such that
\begin{equation*}
(u'',v) + \beta\, (\nabla u',\nabla v) + \alpha\,(u',v) + (\nabla u,\nabla v) = 0
\quad \forall\, v\in H_0^1(\Omega), \ t\in (0,\infty), 
\end{equation*}
with initial conditions
\(
u(0)=u_0,\: u'(0)=u_1 \; \text{in} \; \Omega.
\)

We now state a higher-order regularity and decay result. the operator $A=-\Delta$, we consider the linear self-adjoint positive definite operator on $L^{2}(\Omega)$ having dense domain \( D(A)=H^{2}(\Omega)\cap H_{0}^{1}(\Omega). \) For any $k\ge0$, we define \( \dot{H}^{k}(\Omega):=D\!\left(A^{\frac{k}{2}}\right), \) which is a subspace of $H^{k}(\Omega)$. The corresponding norm is given by \( |\chi|_{k} = \left\| A^{\frac{k}{2}}\chi \right\|. \) The following Poincar\'e inequality will be used throughout the analysis:
\begin{equation*}
\|u\|^{2}\leq \frac{1}{\lambda_{1}}
\left\|A^{\frac{1}{2}}u\right\|^{2},
\end{equation*}
where \(\lambda_{1}\) denotes the smallest positive eigenvalue of the Laplace operator \(A\).

The proof of Theorem \ref{regularity} is established in 
\cite{Krishan2025}.

\begin{theorem}\label{regularity}
Let $j>1$, and assume that
$u(0)\in D(A^j)$ and
$u'(0)\in D(A^{j-\frac12})$.
Define
\[
E_{A^{(j)}}^{(k)}(u)(t)
:=\frac12\left(
|u^{(k)}(t)|_{2j-1}^2
+|u^{(k-1)}(t)|_{2j}^2
\right),
\quad |v|_s:=\|A^{s/2}v\|.
\]
For
\(
0<\delta\leq\frac{\alpha+\beta\lambda_1}{2},
\qquad
\delta<\frac{\lambda_1}{\alpha+\beta\lambda_1},
\)
the estimate
\[
E_{A^{(j)}}^{(k)}(u)(t)
\leq Ce^{-2\delta t}E_{A^{(j)}}^{(k)}(u)(0)
\]
holds for all $t\geq0$.
At $\delta=\lambda_1/(\alpha+\beta\lambda_1)$,
provided $\delta\leq(\alpha+\beta\lambda_1)/2$,
the corresponding bound is
\[
|u^{(k)}(t)|_{2j-1}^2
\leq Ce^{-2\delta t}E_{A^{(j)}}^{(k)}(u)(0).
\]
Here $C$ is independent of $t$.
\end{theorem}

\section{Semidiscrete DG formulation}\label{s3}
In this section, we introduce the finite-dimensional semidiscrete formulation within a discontinuous Galerkin framework. 

\subsection{DG Formulation} We assume a shape-regular family of triangulations $\{\mathcal{T}_h\}_h$ of $\overline{\Omega}$, in the sense of \cite{ciarlet2002finite}, where each triangulation $\mathcal{T}_h$ consists of triangular elements $T$ with diameter $h_T$, and set
$
h:=\max_{T\in\mathcal{T}_h} h_T.
$
Let $\Gamma_h$ denote the collection of all edges of $\mathcal{T}_h$, and let $\bn_e$ be a unit normal vector associated with an edge $e\in\Gamma_h$. For an interior edge $e$ shared by two elements $T_r$ and $T_s$ (with $r<s$), the normal vector $\bn_e$ is chosen to point from $T_r$ to $T_s$. If $e\subset \partial \Omega$, then $\bn_e$ denotes
the outward unit normal to $\partial \Omega$. For a scalar function $\rho$, the average $\{\rho\}$ and jump $[\rho]$ across an interior
edge $e=T_r\cap T_s$ are defined by
\[
\{\rho\}=\frac{1}{2}\big(\rho|_{T_r}+\rho|_{T_s}\big)\big|_e,
\qquad
[\rho]=\big(\rho|_{T_r}-\rho|_{T_s}\big)\big|_e.
\]
On a boundary edge $e\subset\partial\Omega$, both the average and jump are taken as the trace of $\rho$ on $e$. We now introduce the following discontinuous function space:
\begin{equation*}
V := \big\{\phi\in L^2(\Omega):\ \phi|_T\in W^{2,4/3}(T)\ \text{for all } T\in\mathcal{T}_h \big\}.
\end{equation*}
This space is equipped with the broken norm
\begin{equation*}
\vertiii{\phi}:=
\bigg(
\sum_{T\in\mathcal{T}_h}\|\nabla\phi\|_{L^2(T)}^2
+\sum_{e\in\Gamma_h}\frac{\sigma_e}{|e|}\int_e[\phi]^2\,ds
\bigg)^{1/2},
\qquad \forall\,\phi\in V,
\end{equation*}
where $|e|$ denotes the length of the edge $e$, and $\sigma_e>0$ is a penalty parameter
associated with $e\in\Gamma_h$.
We shall frequently employ the following $L^2$-estimate (see \cite[Lemma~6.2]{girault2005discontinuous}): For any
$\phi\in V$
\begin{equation}\label{L2}
\|\phi\|^2_{L^2(\Omega)} \le \gamma\,\vertiii{\phi}^2,
\end{equation}
where  $\gamma>0$ is a constant.

\noindent
Now, define the discrete bilinear form $a_{dG}:\,V\times V\to\mathbb{R}$ by
\begin{eqnarray*}
a_{dG}(\phi,\psi) &:=&
\sum_{T\in\mathcal{T}_h}\int_T  \nabla\phi\cdot \nabla\psi\,d\bx -\sum_{e\in\Gamma_h}\int_e \{ \nabla\phi\}\cdot\bn_e\,[\psi]\,ds
\nonumber \\ 
&& \quad - \sum_{e\in\Gamma_h}\int_e \{ \nabla\psi\}\cdot\bn_e\,[\phi]\,ds  +\sum_{e\in\Gamma_h}\frac{\sigma_e}{|e|}\int_e [\phi][\psi]\,ds,
\quad \forall \; \phi,\,\psi\in V
\end{eqnarray*}
The  symmetric DG  formulation is to seek 
$
u:[0,\infty)\to V
$
such that
\begin{equation}\label{1.5}
(u'',\,v)
+ \beta\,a_{dG}(u',\,v)
+ \alpha\,(u',\,v)
+ a_{dG}(u,\,v)
= 0
\quad \forall\; v \in V, 
\end{equation}
with initial conditions
\(
u(0)=u_0,\: u'(0)=u_1 \; \text{in} \; \Omega.
\)

Next, we define the finite-dimensional subspace $V_h$ of the discontinuous Sobolev space $V$ by
\begin{equation*}
V_h := \big\{\psi\in L^2(\Omega):\ \psi|_T \in \mathbb{P}_r(T)\ \text{for all } T\in\mathcal{T}_h \big\},
\end{equation*}
where $\mathbb{P}_r(T)$ denotes the space of polynomials of degree at most $r$ on the element $T$, with $r\ge 1$.

We next introduce a local $L^2$-projection operator $P_h:H^s(T)\subset L^2(T)\to \mathbb{P}_r(T)$ for $s\in\{1,2,\ldots,r+1\}$. This operator satisfies the standard approximation properties stated below
(see \cite[Section~6.1.5]{riviere2008discontinuous}).
\begin{lemma}\label{prelim}
Let $T\in\mathcal{T}_h$ and $w\in H^s(T)$ with $s\in\{1,2,\ldots,r+1\}$.
Then there exists a local $L^2$-projection $P_h w\in\mathbb{P}_r(T)$ and a constant $C>0$,
independent of $w$ and $h_T$, such that
\begin{equation*}
\|w-P_h w\|_{m,T} \le C\,h_T^{\,s-m}\,|w|_{s,T},
\quad \forall \; m \; \in\{0,1,\ldots,s\}. 
\end{equation*}
\end{lemma}
\noindent

Now, define the discrete bilinear form $a_{dG}:\,V_h\times V_h\to\mathbb{R}$ by
\begin{eqnarray*}
a_{dG}(\phi_h,\psi_h) &:=&
\sum_{T\in\mathcal{T}_h}\int_T  \nabla\phi_h\cdot \nabla\psi_h\,d\bx -\sum_{e\in\Gamma_h}\int_e \{ \nabla\phi_h\}\cdot\bn_e\,[\psi_h]\,ds
\nonumber \\ 
&& \quad - \sum_{e\in\Gamma_h}\int_e \{ \nabla\psi_h\}\cdot\bn_e\,[\phi_h]\,ds  +\sum_{e\in\Gamma_h}\frac{\sigma_e}{|e|}\int_e [\phi_h][\psi_h]\,ds,
\quad \forall \; \phi_h,\,\psi_h\in V_h
\end{eqnarray*}
with the  discrete energy norm
\begin{equation*}
\vertiii{\phi_h}:= \bigg(
\sum_{T\in\mathcal{T}_h}\|\nabla\phi_h\|_{L^2(T)}^2 +\sum_{e\in\Gamma_h}\frac{\sigma_e}{|e|}\int_e[\phi_h]^2\,ds
\bigg)^{1/2}, \quad \forall\; \phi_h \, \in V_h.
\end{equation*}
The semi-discrete symmetric DG  formulation is to seek 
$
u_h:[0,\infty)\to V_h
$
such that
\begin{equation}\label{semiw1}
(u_{h}'',\,v_h)
+ \beta\,a_{dG}(u_{h}',\,v_h)
+ \alpha\,(u_{h}',\,v_h)
+ a_{dG}(u_h,\,v_h)
= 0
\quad \forall\; v_h \in V_h, 
\end{equation}
with the initial conditions
\(
u_h(0)=u_{0h}
\)
and
\(
u_{h}'(0)=u_{1h},
\)
which will be specified later.
By setting
\(
u_{h,\delta}= e^{\delta t} u_h,
\)
we obtain\\
\begin{eqnarray}\label{delta_derivative}
u_{h,\delta}'
- \delta u_{h,\delta} = e^{\delta t} u_{h}',\quad
u_{h,\delta}''-2\delta u_{h,\delta}'+\delta^2 u_{h,\delta}=e^{\delta t}u_h''.
\end{eqnarray}
The semidiscrete problem \eqref{semiw1} can be rewritten as
\begin{eqnarray}\label{semiw2}
(u_{h,\delta}'',v_h)
+ \beta\,a_{dG}(u_{h,\delta}',v_h)
+ (\alpha-2\delta)\,(u_{h,\delta}',v_h)
&+& (1-\beta\delta)\,a_{dG}(u_{h,\delta},v_h)\nonumber\\
&-& \delta(\alpha-\delta)\,(u_{h,\delta},v_h)
= 0,
\quad \forall\; v_h\in V_h.
\end{eqnarray}

We now state the continuity and coercivity properties of the bilinear form
$a_{dG}(\cdot,\cdot)$ on the finite-dimensional space $V_h$. The proofs can be found in \cite[Lemma~3]{wheeler1978elliptic} and \cite[p.~41]{riviere2008discontinuous}.

\begin{lemma}\label{discoer}
For all $\phi_h,\psi_h\in V_h$, the bilinear form $a_{dG}(\cdot,\cdot)$ satisfies
\begin{equation}\label{conteq}
|a_{dG}(\phi_h,\psi_h)| \le \mu\,\vertiii{\phi_h}\,\vertiii{\psi_h}.
\end{equation}
Moreover, if $\sigma_e \ge \sigma_0$ for a sufficiently large constant $\sigma_0>0$, then
\begin{equation}\label{coereq}
a_{dG}(w_h,w_h) \ge \nu\,\vertiii{w_h}^2,
\quad \forall\; w_h\in V_h,
\end{equation}
where the positive constants $\mu$ and $\nu$ are independent of $h$.
\end{lemma}

\noindent
The discrete energy is defined by
\begin{equation*}
E_h(u_h)(t)
= \frac{1}{2}\big(\|u_{h}'\|^2 +\vertiii{u_h(t)}^2\big).
\end{equation*}
Finally, we recall several trace and inverse inequalities for functions in $V_h$
(see \cite[Section~1.4.3]{di2011mathematical}).

Below, we discuss uniform exponential decay property of the semidiscrete solution
\begin{theorem}\label{semidecay}
Let $u_h$ be the solution of \eqref{semiw1}, and further, let
\begin{equation} \label{delta-range}
0<\delta\leq
\min\left\{
\frac{\alpha\gamma+\beta\nu}{2\gamma},
\frac{\nu}{\beta\nu+\alpha\gamma}
\right\}.
\end{equation}
Then the following estimates hold for every $t\geq0$:
\begin{itemize}
\item[(i)]
If
\[
\delta<\frac{\nu}{\beta\nu+\alpha\gamma},
\]
then
\[
E_h(u_h)(t)
\leq C e^{-2\delta t}E_h(u_h)(0).
\]
In particular, this conclusion also holds for
$\delta=(\alpha\gamma+\beta\nu)/(2\gamma)$,
provided that this value is strictly smaller than
$\nu/(\beta\nu+\alpha\gamma)$.

\item[(ii)]
If
\[
\delta=\frac{\nu}{\beta\nu+\alpha\gamma},
\]
then
\[
\|u_h'(t)\|^2
\leq C e^{-2\delta t}E_h(u_h)(0).
\]
\end{itemize}
Here $C>0$ is independent of $h$ and $t$, but may depend on
$\alpha,\beta,\gamma,\mu,\nu$, and $\delta$.
\end{theorem}

\begin{proof}
The argument is based on an energy bound for the exponentially weighted solution $u_{h,\delta}=e^{\delta t}u_h$. Now, introduce
\[
\mathcal{F}_{h,\delta}(t)
:=
\frac12\left(
\|u_{h,\delta}'\|^2
+(1-\beta\delta)a_{dG}(u_{h,\delta},u_{h,\delta})
-\delta(\alpha-\delta)\|u_{h,\delta}\|^2
\right).
\]
Testing \eqref{semiw2} with $u_{h,\delta}'$ shows
\[
\frac{d}{dt}\mathcal{F}_{h,\delta}(t)
+\beta a_{dG}(u_{h,\delta}',u_{h,\delta}')
+(\alpha-2\delta)\|u_{h,\delta}'\|^2
=0.
\]
By  \eqref{L2} and \eqref{coereq}, there holds
\[
a_{dG}(u_{h,\delta}',u_{h,\delta}')
\geq \nu\vertiii{u_{h,\delta}'}^2
\geq \frac{\nu}{\gamma}\|u_{h,\delta}'\|^2.
\]
Consequently,
\begin{equation*}
\frac{d}{dt}\mathcal{F}_{h,\delta}(t)
+\left(\alpha-2\delta+\frac{\beta\nu}{\gamma}\right)
\|u_{h,\delta}'\|^2
\leq0.
\end{equation*}
The coefficient of the last term is nonnegative under the
assumed restriction on $\delta$. Hence
$\mathcal{F}_{h,\delta}(t)\leq\mathcal{F}_{h,\delta}(0)$.\\
To relate this bound to the discrete energy, define
\[
\widetilde{E}_h(u_{h,\delta})(t)
:=
\frac12\left(
\|u_{h,\delta}'\|^2
+(1-\beta\delta)\nu\vertiii{u_{h,\delta}}^2
-\delta(\alpha-\delta)\|u_{h,\delta}\|^2
\right).
\]
Since
$\delta\leq\nu/(\beta\nu+\alpha\gamma)$,
we obtain $1-\beta\delta\geq0$ and a use of the  coercivity therefore yields
\[
\widetilde{E}_h(u_{h,\delta})(t)
\leq \mathcal{F}_{h,\delta}(t).
\]
At the initial time,
\[
u_{h,\delta}(0)=u_h(0),
\qquad
u_{h,\delta}'(0)=u_h'(0)+\delta u_h(0).
\]
Using these identities together with 
 \eqref{L2} and \eqref{conteq}, we arrive at 
\begin{eqnarray}\label{semiw4}
\widetilde{E}_h(u_{h,\delta})(t)
&\leq& \mathcal{F}_{h,\delta}(0) \nonumber\\
&\leq& \frac12\Big(
\|u_h'(0)+\delta u_h(0)\|^2
+(1-\beta\delta)\mu\vertiii{u_h(0)}^2
-\delta(\alpha-\delta)\|u_h(0)\|^2
\Big) \nonumber\\
&\leq& C E_h(u_h)(0).
\end{eqnarray}
Set
\[
\gamma_0
:=\nu-\delta(\beta\nu+\alpha\gamma).
\]
An application of \eqref{L2} to the negative term
$-\alpha\delta\|u_{h,\delta}\|^2$ shows
\[
\widetilde{E}_h(u_{h,\delta})(t)
\geq
\frac12\left(
\|u_{h,\delta}'\|^2
+\delta^2\|u_{h,\delta}\|^2
+\gamma_0\vertiii{u_{h,\delta}}^2
\right).
\]
Moreover,
\[
e^{\delta t}u_h'
=u_{h,\delta}'-\delta u_{h,\delta},
\]
so that
\[
e^{2\delta t}\|u_h'\|^2
\leq
2\left(
\|u_{h,\delta}'\|^2
+\delta^2\|u_{h,\delta}\|^2
\right).
\]
For the part \textup{(i)}, the strict restriction on $\delta$
ensures that $\gamma_0>0$. Combining the preceding bounds,
we now obtain
\begin{eqnarray*}
\widetilde{E}_h(u_{h,\delta})(t)
\geq
\frac{e^{2\delta t}}{2}
\left(
\frac12\|u_h'(t)\|^2
+\gamma_0\vertiii{u_h(t)}^2
\right)
\geq
\min\left\{\frac12,\gamma_0\right\}
e^{2\delta t}E_h(u_h)(t).
\end{eqnarray*}
The assertion now follows from \eqref{semiw4}.

For the part \textup{(ii)}, the endpoint choice of $\delta$
gives $\gamma_0=0$. Although the preceding lower bound
no longer controls the full discrete energy, it still
controls the velocity
\begin{eqnarray}\label{semiw5}
\widetilde{E}_h(u_{h,\delta})(t)
\geq
\frac12\left(
\|u_{h,\delta}'\|^2
+\delta^2\|u_{h,\delta}\|^2
\right)
\geq
\frac14\|u_{h,\delta}'-\delta u_{h,\delta}\|^2
=\frac14 e^{2\delta t}\|u_h'(t)\|^2.
\end{eqnarray}
Combining \eqref{semiw5} with \eqref{semiw4}
proves the second estimate. This concludes the rest of the proof.
\end{proof}

We now extend the decay analysis to higher-order energies.
For an integer $k\geq1$, differentiating \eqref{semiw1}
$k-1$ times with respect to time implies
\begin{equation}\label{semiw6}
\bigl(u_h^{(k+1)},\phi_h\bigr)
+ \beta\, a_{dG}\bigl(u_h^{(k)},\phi_h\bigr)
+ \alpha \bigl(u_h^{(k)},\phi_h\bigr)
+ a_{dG}\bigl(u_h^{(k-1)},\phi_h\bigr)=0,
\quad \forall \; \phi_h\in V_h.
\end{equation}
\noindent
Set $z_{h,\delta}=e^{\delta t}u_h^{(k-1)}$. Since
\begin{eqnarray*}
e^{\delta t}u_h^{(k)}
=z_{h,\delta}'-\delta z_{h,\delta}, \quad 
e^{\delta t}u_h^{(k+1)}
=z_{h,\delta}''-2\delta z_{h,\delta}'
+\delta^2z_{h,\delta},
\end{eqnarray*}
the equation \eqref{semiw6} becomes
\begin{eqnarray}\label{semiw7}
\bigl(z_{h,\delta}'',\phi_h\bigr)
&+&\beta a_{dG}\bigl(z_{h,\delta}',\phi_h\bigr)
+(\alpha-2\delta)\bigl(z_{h,\delta}',\phi_h\bigr)
\nonumber\\
&+&(1-\beta\delta)a_{dG}\bigl(z_{h,\delta},\phi_h\bigr)
-\delta(\alpha-\delta)
\bigl(z_{h,\delta},\phi_h\bigr)=0,
\quad \phi_h\in V_h.
\end{eqnarray}

By \eqref{semiw7}, $z_h$ satisfies
the same semidiscrete equation  \eqref{semiw2} by replacing  $u_h$ by $z_h$, therefore, the conclusions of Theorem ~\ref{semidecay} are valid. Thus, we obtain the following theorem.
\begin{theorem}
Let $u_h$ be a solution of the semidiscrete problem \eqref{semiw1}, and let $k\geq1$ be an integer
for which the required time derivatives exist. Define
\begin{equation*}
E^{(k)}(u_h)(t)
:=\frac12\|u_h^{(k)}(t)\|^2
+\vertiii{u_h^{(k-1)}(t)}^2.
\end{equation*}
Assume that \eqref{delta-range} holds.
Then, for all $t\geq0$
\begin{itemize}
\item[(i)] If $\delta<\nu/(\beta\nu+\alpha\gamma)$, then
\begin{equation*}
E^{(k)}(u_h)(t)
\leq C e^{-2\delta t}E^{(k)}(u_h)(0).
\end{equation*}
\item[(ii)] If $\delta=\nu/(\beta\nu+\alpha\gamma)$, then
\begin{equation*}
\|u_h^{(k)}(t)\|^2
\leq C e^{-2\delta t}E^{(k)}(u_h)(0).
\end{equation*}
\end{itemize}
Here, the constant $C$ is independent of $h$ and $t$, but may
depend on $\alpha,\beta,\gamma,\mu,\nu$, and $\delta$.
\end{theorem}

\subsection{Error Estimates}
In this subsection, we establish error bounds for the semidiscrete approximation and investigate their long-time behavior. To this end, let \(R_h u\in V_h\) be the elliptic projection of the exact solution \(u\), characterized by
\begin{eqnarray}\label{elpro1}
a_{dG}(u-R_h u,\,\phi_h)=0
\qquad \forall\,\phi_h\in V_h.
\end{eqnarray}
Differentiating \eqref{elpro1} with respect to time, we have
\begin{eqnarray*}
a_{dG}((u-R_h u)',\,\phi_h)&=&0
\qquad \forall\,\phi_h\in V_h,
\\
a_{dG}((u-R_h u)'',\,\phi_h)&=&0
\qquad \forall\,\phi_h\in V_h,
\end{eqnarray*}
where $R_h u(0)=P_h u_0$ and $(R_h u)'(0)=P_h u_1$.
To facilitate the analysis, we split the error between the exact and semidiscrete solutions into two components:
\[
e:=u-u_h=(u-R_h u)+(R_h u-u_h)=\xi+\rho,
\]
where \(\xi=u-R_h u\)  and  \(\rho=R_h u-u_h\).
 
The defining property of the elliptic projection immediately implies that
\begin{equation*}
a_{dG}(\xi_\delta', \phi_h) = \delta e^{\delta t} a_{dG}(\xi, \phi_h)
+ e^{\delta t} a_{dG}((u' - R_h u'), \phi_h) = 0, \quad \forall \; \phi_h \in V_h .
\end{equation*}
Applying Lemma \ref{prelim} and following the techniques utilized in \cite[Lemma 3.8]{BajpaiSwain2025}, we have the following results for the projection $R_h$
\begin{eqnarray}
\sum_{m=0}^{2}
\left(\|\xi^{(m)}\|+h\vertiii{\xi^{(m)}}\right)
&\leq& Ch^{r+1}
\sum_{m=0}^{2}\|u^{(m)}\|_{r+1},\label{approx1}\\
\sum_{m=0}^{2}
\left(\|\xi_\delta^{(m)}\|
+h\vertiii{\xi_\delta^{(m)}}\right)
&\leq& C h^{r+1}e^{\delta t}
\sum_{m=0}^{2}\|u^{(m)}\|_{r+1}\label{approx2}.
\end{eqnarray}
Here $\xi^{(0)}=\xi$, $\xi^{(1)}=\xi'$,
$\xi^{(2)}=\xi''$, and similarly for $\xi_\delta=e^{\delta t}\xi$.
The second estimate follows from
\begin{equation*}
\xi_\delta^{(m)}
=e^{\delta t}\sum_{j=0}^{m}
\binom{m}{j}\delta^{m-j}\xi^{(j)},
\qquad m=0,1,2.
\end{equation*}
The constant $C>0$ is independent of $h$ and $t$.\\
The error equation follows by subtracting \eqref{semiw1}
from \eqref{1.5} and applying the orthogonality property
of the elliptic projection:
\begin{equation}\label{semiw07}
(\rho'', \phi_h) + \beta\, a_{dG}(\rho', \phi_h) + \alpha(\rho', \phi_h)
+ a_{dG}(\rho, \phi_h)
= -(\xi''+\alpha\xi', \phi_h),
\quad \forall\; \phi_h \in V_h .
\end{equation}
Setting $\rho_\delta=e^{\delta t}\rho$ and multiplying
\eqref{semiw07} by $e^{\delta t},$  we arrive at
\begin{eqnarray}\label{semiw8}
(\rho_\delta'', \phi_h) &+& \beta\, a_{dG}(\rho_\delta', \phi_h)
+ (\alpha - 2\delta)(\rho_\delta', \phi_h)
+ (1 - \beta\delta)\, a_{dG}(\rho_\delta, \phi_h)\nonumber\\
& -& \delta(\alpha - \delta)(\rho_\delta, \phi_h)
= -e^{\delta t}(\xi'', \phi_h) - \alpha e^{\delta t}(\xi', \phi_h), \quad \forall\;  \phi_h \in V_h .
\end{eqnarray}
\begin{theorem}\label{semierr1}
Let $\rho$ satisfy \eqref{semiw07}, and choose
$0<\delta_0<\alpha\gamma+\beta\nu$. Further, let
\begin{equation*}
0<\delta\leq
\min\left\{
\frac{\alpha\gamma+\beta\nu-\delta_0}{2\gamma},
\frac{\nu}{\beta\nu+\alpha\gamma}
\right\}.
\end{equation*}
Then, for every $t\geq0$, the following estimates hold:
\begin{enumerate}[(i)]
\item If \[\delta<\frac{\nu}{(\beta\nu+\alpha\gamma)},\] then
\begin{equation*}
E_h(\rho)(t)
\le
C e^{-2\delta t} E_h(\rho)(0)
+
C \int_{0}^{t} e^{-2\delta (t-s)}
\left(
\|\xi''(s)\|^{2} + \|\xi'(s)\|^{2}
\right)\, ds.
\end{equation*}
This includes
$\delta=(\alpha\gamma+\beta\nu-\delta_0)/(2\gamma)$
whenever it is strictly below $\nu/(\beta\nu+\alpha\gamma)$.

\item If
\[
0 < \delta = \dfrac{\nu}{\beta \nu+\alpha \gamma},
\]
then
\begin{equation*}
\|\rho'(t)\|^{2}
\le
C e^{-2\delta t} E_h(\rho)(0)
+
C \int_{0}^{t} e^{-2\delta (t-s)}
\left(
\|\xi''(s)\|^{2} + \|\xi'(s)\|^{2}
\right)\, ds.
\end{equation*}
\end{enumerate}
The constant $C$ is independent of $h$ and $t$, but may
depend on $\alpha,\beta,\gamma,\mu,\nu,\delta$, and $\delta_0$.
\end{theorem}
\begin{proof}
For $\rho_\delta=e^{\delta t}\rho$, we introduce
\begin{equation*}
F_{h,\delta}(t)
:=\frac12\left(
\|\rho_\delta'\|^2
+(1-\beta\delta)a_{dG}(\rho_\delta,\rho_\delta)
-\delta(\alpha-\delta)\|\rho_\delta\|^2
\right).
\end{equation*}
Taking $\phi_h=\rho_\delta'$ in \eqref{semiw8} and using
the symmetry of $a_{dG}$ gives
\begin{equation*}
\frac{d}{dt}F_{h,\delta}(t)
+\beta a_{dG}(\rho_\delta',\rho_\delta')
+(\alpha-2\delta)\|\rho_\delta'\|^2
=-e^{\delta t}(\xi''+\alpha\xi',\rho_\delta').
\end{equation*}
By \eqref{coereq} and \eqref{L2},
\begin{equation*}
a_{dG}(\rho_\delta',\rho_\delta')
\geq\frac{\nu}{\gamma}\|\rho_\delta'\|^2.
\end{equation*}
The forcing term is bounded by Young's inequality
\begin{eqnarray*}
e^{\delta t}
\left|(\xi''+\alpha\xi',\rho_\delta')\right|
&\leq&
\frac{\delta_0}{2\gamma}\|\rho_\delta'\|^2
+\frac{\gamma}{2\delta_0}
e^{2\delta t}\|\xi''+\alpha\xi'\|^2\\
&\leq&
\frac{\delta_0}{2\gamma}\|\rho_\delta'\|^2
+\frac{\gamma}{\delta_0}e^{2\delta t}
\left(\|\xi''\|^2+\alpha^2\|\xi'\|^2\right).
\end{eqnarray*}
Since
$\alpha-2\delta+\beta\nu/\gamma\geq\delta_0/\gamma$,
we obtain
\begin{equation*}
\frac{d}{dt}F_{h,\delta}(t)
+\frac{\delta_0}{2\gamma}\|\rho_\delta'\|^2
\leq
\frac{\gamma}{\delta_0}e^{2\delta t}
\left(\|\xi''\|^2+\alpha^2\|\xi'\|^2\right).
\end{equation*}

Employing a kick-back argument  and \eqref{conteq}
followed by integration with respect to time from $0$ to $t$, we arrive at
\begin{eqnarray*}
\widetilde{E}_h(\rho_\delta)(t)
 \le
\frac{1}{2}\Big(
\|\rho_\delta'(0)\|^2
+ (1-\beta\delta)\mu\,\vertiii{\rho_\delta(0)}^2
&-& (\alpha-\delta)\delta\,\|\rho_\delta(0)\|^2
\Big) \\
&  +&
\frac{\gamma e^{\delta t}}{\delta_0} \int_{0}^{t}
e^{2\delta s} \left(\|\xi''(s)\|^{2} 
 +
\alpha^{2}\|\xi'(s)\|^{2}
\right)\, ds.
\end{eqnarray*}
The remainder of the proof proceeds along the same lines as the argument used in the proof of Theorem~\ref{semidecay}.
\end{proof}
\noindent
With $u_{0h}=R_hu_0$, we have $\rho(0)=0$ and hence
\[
E_h(\rho)(0)=\frac12\|\rho'(0)\|^{2}.
\]
Choose $u_{1h}$ as the $L^2$-projection of $u_1$, or as
an interpolant with the same $L^2$ approximation order.
Since $\rho'(0)=R_hu_1-u_{1h}$, the projection and
interpolation estimates give
\[
E_h(\rho)(0)\le C h^{2r+2}\|u_1\|_{r+1}^{2}.
\]
Therefore, applying Theorem~\ref{semierr1} yields
\begin{equation}\label{semiw12}
\|\rho'(t)\|^{2}+\vertiii{\rho(t)}^{2}
\le
C h^{2r+2} e^{-2\delta t}
\left(
\|u_1\|_{r+1}^{2}
+
\int_{0}^{t} e^{2\delta s}
\left(
\|u''(s)\|_{r+1}^{2}
+
\|u'(s)\|_{r+1}^{2}
\right)\, ds
\right).
\end{equation}
Applying Theorem~\ref{regularity},   
Taking $j=\tfrac{r+2}{2}$ and $k=2$, we obtain
\begin{equation}\label{semiw13}
\|u''(t)\|_{r+1}^{2}
\le
C e^{-2\delta t} E_{A^{(\frac{r+2}{2})}}^{(2)}(u)(0) \le C e^{-2\delta t} \big(|u^{(2)}(0)|_{r+1}^{2} +
|u'(0)|_{r+2}^{2} \big).
\end{equation}
Evaluating \eqref{1.1} at $t=0$ gives
\begin{equation}\label{semiw14}
u''(0)=\beta\Delta u_1-\alpha u_1+\Delta u_0.
\end{equation}
Thus
\begin{equation*}
\|u''(0)\|_{r+1}^2
\leq C\left(
\beta^2\|u_1\|_{r+3}^2
+\alpha^2\|u_1\|_{r+1}^2
+\|u_0\|_{r+3}^2
\right),
\end{equation*}
and \eqref{semiw13} implies
\begin{equation}\label{semiw15}
\|u''(t)\|_{r+1}^2
\leq Ce^{-2\delta t}
\left(
\beta^2\|u_1\|_{r+3}^2
+\|u_1\|_{r+2}^2
+\|u_0\|_{r+3}^2
\right).
\end{equation}
For $j=\tfrac{r+2}{2}$ and $k=1$, Theorem \ref{regularity} similarly gives
\begin{eqnarray}\label{semiw16}
\|u'(t)\|_{r+1}^{2}
&\le&
C e^{-2\delta t} E_{A^{(\frac{r+2}{2})}}^{(1)}(u)(0) \nonumber\\
&\le&
C e^{-2\delta t}\big(
\|u_1\|_{r+1}^{2}
+
\|u_0\|_{r+2}^{2}
\big).
\end{eqnarray}
Substituting the bounds \eqref{semiw13}--\eqref{semiw16} into \eqref{semiw12}, we obtain
\begin{equation}\label{maxnorm}
\|\rho'(t)\|^{2}+\vertiii{\rho(t)}^{2}
\le
C h^{2r+2} (1+t)e^{-2\delta t}
\big(
\beta^2 \|u_1\|_{r+3}^{2}
+
\|u_1\|_{r+2}^{2}
+
\|u_0\|_{r+3}^{2}
\big).
\end{equation}
Finally, utilizing the triangle inequality together with the approximation estimates for $\xi$, we conclude
\begin{equation*}
\|u'(t)-u_h'(t)\|^{2}
+
h^2\vertiii{u(t)-u_h(t)}^{2}
\le
C h^{2r+2} (1+t)e^{-2\delta t}
\big(
\beta^2 \|u_1\|_{r+3}^{2}
+
\|u_1\|_{r+2}^{2}
+
\|u_0\|_{r+3}^{2}
\big).
\end{equation*}

\begin{remark}
If $\beta\to0$, then the above estimate reduces accordingly, leading to
\begin{equation*}
\|u'(t)-u_h'(t)\|^{2}
+
h^2\vertiii{u(t)-u_h(t)}^{2}
\le
C h^{2r+2} (1+t)e^{-2\delta t}
\big(
\|u_1\|_{r+2}^{2}
+
\|u_0\|_{r+3}^{2}
\big).
\end{equation*}
\end{remark}
Below, we present an optimal error estimate in the energy norm as a theorem.
\begin{theorem}\label{semierr2}
Let $u_{0h}$ and $u_{1h}$ be chosen as either the $L^2$--projection or the finite element interpolant of $u_0$ and $u_1$, respectively.
Choose
$0<\delta_0<\alpha\gamma+\beta\nu$ and
\(
0<\delta\leq
\frac{\alpha\gamma+\beta\nu-\delta_0}{2\gamma},
\quad
\delta<\frac{\nu}{\beta\nu+\alpha\gamma}.
\)
For $t\ge0$,
\[
\vertiii{(u-u_h)(t)}
\le
C h^{r} (1+t)e^{-\delta t}
\big(
\|u_0\|_{r+2}
+\beta\|u_1\|_{r+2}
+\|u_1\|_{r+1}
\big).
\]
Moreover, if $\beta\to0$ and $\alpha>0$, then
\[
\vertiii{(u-u_h)(t)}
\le
C h^{r} (1+t)e^{-\delta t}
\big(
\|u_0\|_{r+2}
+\|u_1\|_{r+1}
\big).
\]
\end{theorem}
\begin{proof}
Since $u_{0h}$ and $u_{1h}$ are taken as either the $L^2$--projection or the interpolant of $u_0$ and $u_1$, standard approximation results yield
\[
E_h(\rho)(0)
\le
C h^{2r}
\big(
\|u_0\|_{r+1}^{2}
+
\|u_1\|_{r}^{2}
\big).
\]
Consequently, applying Theorem~\ref{semierr1}, we obtain for
\(
0 < \delta < \min\left(\frac{\alpha\gamma + \beta \nu - \delta_{0}}{2\gamma}, \, \dfrac{1}{\beta \nu+\alpha \gamma}\right)
\quad \text{or} \quad
0 < \delta= \frac{\alpha\gamma + \beta \nu - \delta_{0}}{2\gamma},
\)
that
\begin{eqnarray}\label{semiw17}
\vertiii{\rho(t)}^{2}
&\le &
C e^{-2\delta t} E_h(\rho)(0)
+
C  e^{-2\delta t}
\left(
\int_{0}^{t} e^{2\delta s}
\big(
\|\xi''(s)\|^{2}
+
\|\xi'(s)\|^{2}
\big)\, ds
\right)\nonumber\\
&\leq & C h^{2r} e^{-2\delta t}
\left(\|u_0\|_{r+1}^{2}+
\|u_1\|_{r}^{2}
+
\int_{0}^{t} e^{2\delta s}
\left(
\|u''(s)\|_{r}^{2}
+
\|u'(s)\|_{r}^{2}
\right)\, ds
\right).
\end{eqnarray}
To estimate the integral terms, we apply Theorem~\ref{regularity} to obtain  
\begin{equation}\label{semiw18}
\|u''(s)\|_{r}^{2}
\le
C e^{-2\delta s} E_A^{(2)}(u)(0)
\le
C e^{-2\delta s}
\big(
\beta^2\|u_1\|_{r+2}^{2}
+
\|u_1\|_{r+1}^{2}
+
\|u_0\|_{r+2}^{2}
\big).
\end{equation}
Similarly, 
\begin{equation}\label{semiw19}
\|u'(s)\|_{r}^{2}\le
C e^{-2\delta s} E_A^{(1)}(u)(0)
\le
C e^{-2\delta s}
\big(
\|u_1\|_{r}^{2}
+
\|u_0\|_{r+1}^{2}
\big).
\end{equation}
Substituting \eqref{semiw18} and \eqref{semiw19} into \eqref{semiw17} completes the proof.
\end{proof}

As a consequence of Theorem ~\ref{semierr1} with $u_{0h}=R_h u_0$, one obtains optimal error estimate in $L^2$ norm, but with higher regularity results on the exact solution. Below, we provide an optimal estimate  with or without using elliptic projection as our initial approximation $u_{0h}$ and also using less regularity.

\begin{theorem}\label{semierr3}
Let $u$ and $u_h$ denote the solutions of \eqref{1.1} and \eqref{semiw1}, respectively. Further, let $u_{0h}$ and $u_{1h}$ be chosen as either the $L^2$--projection or the finite element interpolant of $u_0$ and $u_1$, respectively.
Then, under the assumptions on $\delta_0$ and $\delta$ as in theorem ~\ref{semierr2},  there exists a positive constant $C$, independent of $h$, such that
\[
\|u(t)-u_h(t)\|
\le
C h^{r+1}\, t^{1/2} e^{-\delta t}
\big(
\|u_0\|_{r+2}+\|u_1\|_{r+1}
\big).
\]
\end{theorem}
\begin{proof}
Integrating \eqref{semiw07} with respect to time, we obtain
\begin{eqnarray}\label{semiw20}
(\rho',\phi_h)
+ \beta a_{dG}(\rho,\phi_h)
+ \alpha(\rho,\phi_h)
+ a_{dG}(\hat{\rho},\phi_h)
&=&
-(\xi'+\alpha\xi,\phi_h)
\nonumber\\
&
+&(\alpha e(0)+e'(0),\phi_h)
+\beta a_{dG}(\rho(0),\phi_h).
\end{eqnarray}
We choose $u_{0h}$ and $u_{1h}$ as the elliptic projection and the $L^2$--projection of
$u_0$ and $u_1$, respectively. With this choice, $(e'(0),\phi_h)=0$ and $\rho(0)=0$.\\

\noindent
Let $v=\rho+\tilde{\delta}\hat{\rho}$ in \eqref{semiw20}. Using the definition of the
energy  together with Young’s inequality, \eqref{L2},  \eqref{conteq}, and \eqref{coereq}, we obtain
\begin{equation}\label{semiw21}
\frac{d}{dt}\widetilde{E}_{\tilde{\delta}}(\rho)(t)
+
F(\rho)(t)
\le
C\big(
\|\xi'(t)\|^{2}
+\|\xi(t)\|^{2}
+\|\xi(0)\|^{2}
\big),
\end{equation}
where
\[
\widetilde{E}_{\tilde{\delta}}(\rho)(t)
=
\frac12\Big(
\|\rho\|^{2}
+
a_{dG}(\hat{\rho}, \hat{\rho})
\Big)
+
\tilde{\delta}(\rho,\hat{\rho})
+
\frac{\beta\tilde{\delta}}{2}a_{dG}(\hat{\rho},\, \hat{\rho})
+
\frac{\alpha\tilde{\delta}}{2}\|\hat{\rho}\|^{2},
\]
with
\[
E(\rho)(t)
=
\frac12\Big(
\|\rho(t)\|^{2}
+
\vertiii{\hat{\rho}(t)}^{2}
\Big),
\]
and
\[
F(\rho)(t)
=
(\tilde{\delta}-\delta_0)a_{dG}(\hat{\rho},\,\hat{\rho})
+
(\alpha+\frac{\beta \nu}{\gamma}-\delta_0-\tilde{\delta})\|\rho\|^{2}.
\]
By  application of the Young’s inequality, with
\(
0<\tilde{\delta}
\le
\min\left(
\frac{(\alpha\gamma+\beta\nu)}{2\gamma},
\,
\frac{\nu}{4(\alpha\gamma+\beta\nu)}\right)
\), \eqref{L2} and \eqref{coereq}, we arrive at
\begin{eqnarray*}
\widetilde{E}_{\tilde{\delta}}(\rho)(t)
&=&
\frac12\Big(
\|\rho\|^{2}
+
a_{dG}(\hat{\rho}, \hat{\rho})
\Big)
+
\tilde{\delta}(\rho,\hat{\rho})
+
\frac{\beta\tilde{\delta}}{2}a_{dG}(\hat{\rho},\, \hat{\rho})
+
\frac{\alpha\tilde{\delta}}{2}\|\hat{\rho}\|^{2}
\nonumber\\
& \leq &\frac12\Big(
\|\rho\|^{2}
+
a_{dG}(\hat{\rho}, \hat{\rho})
\Big)
+
\frac{\gamma\widetilde\delta}
     {2(\alpha\gamma+\beta\nu)}\|\rho\|^2
+\frac{(\alpha\gamma+\beta\nu)\widetilde\delta}
      {2\nu}
a_{dG}(\widehat\rho,\widehat\rho)
+
\left(\frac{\beta\tilde{\delta}}{2}
+\frac{\alpha\tilde{\delta}}{2}\frac{\gamma}{\nu}\right)
a_{dG}(\widehat\rho,\widehat\rho).\nonumber\\
& \leq &\frac34\Big(
\|\rho\|^{2}
+
a_{dG}(\hat{\rho}, \hat{\rho})
\Big)+\left(\frac{\gamma\tilde{\delta}}{2(\alpha\gamma+\beta\nu)}-\frac14\right)\|\rho\|^2
+\left(\frac{(\alpha\gamma+\beta\nu)\tilde{\delta}}{\nu}-\frac{1}{4}\right)a_{dG}(\widehat\rho,\widehat\rho)\nonumber\\
&\leq &\frac34\Big(
\|\rho\|^{2}
+
a_{dG}(\widehat\rho,\widehat\rho)
\Big).
\end{eqnarray*}
Applying Young's inequality, \eqref{L2} and \eqref{coereq} with $0<\tilde{\delta}\leq \min\left(
\frac{(\alpha\gamma+\beta\nu)}{2\gamma},
\,
\frac{\nu}{2(\alpha\gamma+\beta\nu)} \right)$, we find that
\begin{eqnarray*}
\widetilde{E}_{\tilde{\delta}}(\rho)(t)
&= &
\frac12\Big(
\|\rho\|^{2}
+
a_{dG}(\hat{\rho}, \hat{\rho})
\Big)
+
\tilde{\delta}(\rho,\hat{\rho})
+
\frac{\beta\tilde{\delta}}{2}a_{dG}(\hat{\rho},\, \hat{\rho})
+
\frac{\alpha\tilde{\delta}}{2}\|\hat{\rho}\|^{2}\nonumber\\
&\ge &\frac12\Big(
\|\rho\|^{2}
+a_{dG}(\widehat\rho,\widehat\rho)\Big)
-\left(\frac{\gamma\tilde{\delta}}{2(\alpha\gamma+\beta\nu)}-\frac14\right)\|\rho\|^2
-\left(\frac{(\alpha\gamma+\beta\nu)\tilde{\delta}}{2\nu}-\frac{1}{4}\right)a_{dG}(\widehat\rho,\widehat\rho)
\nonumber\\
&\ge &\frac14\Big(
\|\rho\|^{2}
+a_{dG}(\widehat\rho,\widehat\rho)\Big).
\end{eqnarray*}
Hence,
\begin{eqnarray}\label{semiw22}
\frac14\min\left(1,\nu\right)\Big(
\|\rho\|^{2} +\vertiii{\hat{\rho}}^2\Big) &\leq &\frac14\Big(\|\rho\|^{2} + a_{dG}(\widehat\rho,\widehat\rho)\Big) \nonumber \\
& \le &\widetilde{E}_{\tilde{\delta}}(\rho)(t)
\leq \frac34\Big(\|\rho\|^{2} + a_{dG}(\widehat\rho,\widehat\rho) \Big) \nonumber \\
& \leq &\frac34\max\left(1,\mu\right)\Big(\|\rho\|^{2}
+ \vertiii{\hat{\rho}}^2\Big)
\end{eqnarray}
provided $0<\tilde{\delta}\leq  \min\left(\frac{(\alpha\gamma+\beta\nu)}{2\gamma},
\,
\frac{\nu}{4((\alpha\gamma+\beta\nu)} \right)$.\\

\noindent
Furthermore, from the definition of $F(\rho)(t)$ and for
$0<\delta_0<\tilde{\delta}\le\frac{\alpha\gamma+\beta \nu}{2\gamma}$, $q=\frac{\min\left(1,\nu\right)
}{\max\left(1,\mu\right)}$, it follows that
\begin{eqnarray}\label{semiw23}
F(\rho)(t)
&= &
(\tilde{\delta}-\delta_0 )a_{dG}(\hat{\rho},\,\hat{\rho})
+
\left(\alpha+\frac{\beta \nu}{\gamma}-\delta_0-\tilde{\delta}\right)\|\rho\|^{2}.
\nonumber\\
&\ge &
(\tilde\delta-\delta_0)
\left(\|\rho\|^{2}
+ a_{dG}
(\hat{\rho},\,\hat{\rho})\right)\nonumber\\
&\ge &
(\tilde\delta-\delta_0)\min\left(1,\nu\right)
\left(\|\rho\|^{2}
+ 
\vertiii{\hat{\rho}}^2\right)\nonumber\\
&\ge &
(\tilde\delta-\delta_0)\frac43\frac{\min\left(1,\nu\right)
}{\max\left(1,\mu\right)}\widetilde{E}_{\tilde{\delta}}(\rho)(t)
\nonumber\\
&= &(\tilde\delta-\delta_0)\frac43q\widetilde{E}_{\tilde{\delta}}(\rho)(t)
.
\end{eqnarray}
Substituting \eqref{semiw22} and \eqref{semiw23} into \eqref{semiw21} and integrating in time,
we obtain
\begin{eqnarray}
\widetilde{E}_{\tilde{\delta}}(\rho)(t)
&\le &
e^{-(\tilde\delta-\delta_0)\frac43qt}
\widetilde{E}_{\tilde{\delta}}(\rho)(0)
+
C\int_{0}^{t}
e^{-(\tilde\delta-\delta_0)\frac43q(t-s)}
\big(
\|\xi'(s)\|^{2}
+\|\xi(s)\|^{2}
+\|\xi(0)\|^{2}
\big)\,ds \nonumber\\
&\le &
e^{-2\delta t}
\widetilde{E}_{\tilde{\delta}}(\rho)(0)
+
C\int_{0}^{t}
e^{-2\delta(t-s)}
\big(
\|\xi'(s)\|^{2}
+\|\xi(s)\|^{2}
+\|\xi(0)\|^{2}
\big)\,ds,
\label{semiw24}
\end{eqnarray}
where $\delta=(\tilde\delta-\delta_0)\frac23p$.\\
Since $E(\rho)(0)=0$ and 
applying the estimate \eqref{semiw16}, we have
\[
\|\xi'(s)\|^{2}+\|\xi(s)\|^{2}
\le
C h^{2r+2}
\big(
\|u'(s)\|_{r+1}^{2}
+\|u(s)\|_{r+1}^{2}
\big)
\le
C h^{2r+2} e^{-2\delta s}
\big(
\|u_1\|_{r+1}^{2}
+\|u_0\|_{r+2}^{2}
\big).
\]
Substituting this bound into \eqref{semiw24} yields
\[
\|\rho(t)\|^{2}
\le
C h^{2r+2} t e^{-2\delta t}
\big(
\|u_1\|_{r+1}^{2}
+\|u_0\|_{r+2}^{2}
\big),
\]
and hence
\[
\|u(t)-u_h(t)\|
\le
C h^{r+1} t^{1/2} e^{-\delta t}
\big(
\|u_1\|_{r+1}
+\|u_0\|_{r+2}
\big).
\]
This completes the proof.
\end{proof}
The following theorem provides an error estimate in the maximum norm.
\begin{theorem}\label{semimax}
Let $u_{0h}=R_hu_0$, and let $u_{1h}$ be the
$L^2$-projection or a finite element interpolant of $u_1$.
Under the required regularity  $u_0\in H^{r+4}(\Omega)$ and
$u_1\in H^{r+3}(\Omega)$, the following estimate holds:
\begin{equation*}
\|u(t)-u_h(t)\|_{L^\infty(\Omega)}
\leq Ch^{r+1}(1+|\log h|)e^{-\delta t}
\left(
\|u_0\|_{r+4}+(1+\beta)\|u_1\|_{r+3}
\right).
\end{equation*}
\end{theorem}

\begin{proof}
The discrete Sobolev inequality
\cite[pp.~43]{Brenner2004}, together with \eqref{maxnorm},
gives
\begin{eqnarray*}
\|\rho(t)\|_{L^\infty}
&\leq& C(1+|\log h|)^{1/2}\vertiii{\rho(t)}
\nonumber\\
&\leq& Ch^{r+1}(1+|\log h|)^{1/2}(1+t)^{1/2}e^{-\delta t}
\left(
\beta\|u_1\|_{r+3}
+\|u_1\|_{r+2}
+\|u_0\|_{r+3}
\right).
\end{eqnarray*}
For the projection error, the assumed maximum-norm
estimate, \cite[Corollary~1]{KanschatRannacher2002},
and Theorem \ref{regularity} yield
\begin{eqnarray*}
\|\xi(t)\|_{L^\infty}
&\leq& Ch^{r+1}(1+|\log h|)
\|u(t)\|_{W^{r+1,\infty}}\\
&\leq& Ch^{r+1}(1+|\log h|)e^{-\delta t}
\left(
\|u_1\|_{r+3}+\|u_0\|_{r+4}
\right).
\end{eqnarray*}
For $d=2$, an application of the triangle inequality
then completes the proof.
\end{proof}

\section{Completely Discrete Scheme}\label{s4}

This section introduces a fully discrete approximation obtained by combining a finite-difference discretization in time with a spatially semidiscrete finite element scheme. It is  further shown that the resulting fully discrete solution inherits a uniform exponential decay property analogous to that of the continuous problem.

Let $k>0$ denote the time step size and define the discrete time levels by $t_n = nk$, for $n \ge 0$. For any sufficiently smooth function $\phi$, we denote
$\phi^n := \phi(t_n)$. The backward and forward finite difference operators are defined by
\(
\bar{\partial}_t \phi^n := \frac{\phi^n - \phi^{n-1}}{k}, \;
\partial_t \phi^n := \frac{\phi^{n+1} - \phi^n}{k},
\)
with  $\bar{\partial}_t^0 \phi^n = \phi^n$.
For integers $j \ge 0$, we recursively define higher-order discrete derivatives by
\[
\bar{\partial}_t^{\,j+1} \phi^n
:= \frac{\bar{\partial}_t^{\,j} \phi^n - \bar{\partial}_t^{\,j} \phi^{n-1}}{k}.
\]
We also introduce the midpoint average
\[
\phi^{n+\frac12} := \frac{\phi^{n+1} + \phi^n}{2}.
\]
The second-order central difference operator in time is given by
\[
\partial_t \bar{\partial}_t \phi^n
:= \frac{\phi^{n+1} - 2\phi^n + \phi^{n-1}}{k^2}
= \frac{1}{k}\bigl(\partial_t \phi^n - \bar{\partial}_t \phi^n\bigr).
\]
The fully discrete approximation $U^n \in V_h$ to the exact solution $u(t_n)$ is defined for $n \ge 1$ as
\begin{equation}\label{fully1}
(\partial_t \bar{\partial}_t U^n, \chi_h)
+ \beta\, a_{dG}(\partial_t U^n, \chi_h)
+ \alpha (\partial_t U^n, \chi_h)
+ a_{dG}(U^{n+1}, \chi_h)
= 0,
\qquad \forall \chi_h \in V_h,
\end{equation}
where $U^0=u_{0,h}\in V_h$ be an appropriate approximation
of $u_0$. The value $U^1\in V_h$ is determined by solving
\begin{eqnarray}\label{fullyerrstartup}
\frac{2}{k}(\partial_tU^0,\chi_h)
+\beta a_{dG}(\partial_tU^0,\chi_h)
+\alpha(\partial_tU^0,\chi_h)+a_{dG}(U^1,\chi_h)
=\frac{2}{k}(u_1,\chi_h),
\qquad \chi_h\in V_h,
\end{eqnarray}
where $k$ is the time step and
$\partial_tU^0=(U^1-U^0)/k$.
We associate with the discrete solution the energy functional
\[
E^n(U)
:= \frac12 \left( \|\partial_t U^n\|^2 + \vertiii{U^{n+1}}^2 \right).
\]
The following lemma establishes fundamental properties of the discrete energy and plays
a crucial role in proving the stability and decay behavior of the fully discrete scheme.
\begin{lemma}\label{fullyD1}
Let $U^{n+1}$, for $n\geq1$, satisfy \eqref{fully1}.
Then, for every $0\leq\vartheta\leq1$,
\begin{eqnarray}\label{4.2}
\frac12\bar{\partial}_t
\left(
\|\partial_tU^n\|^2+a_{dG}(U^{n+1},U^{n+1})
\right)
&\leq &
-\vartheta\left(\frac{\alpha\gamma+\beta\nu}{\gamma}\right)
\|\partial_tU^n\|^2
-(1-\vartheta)\alpha\|\partial_tU^n\|^2
\nonumber\\
&-&(1-\vartheta)\beta\nu\vertiii{\partial_tU^n}^2
-\frac{k}{2}\left(
\|\bar{\partial}_t\partial_tU^n\|^2
+\nu\vertiii{\partial_tU^n}^2
\right).
\end{eqnarray}
Moreover, for a modified energy
\begin{equation}\label{4.3}
\widetilde{E}^n_\delta(U)
=\frac12\left(
\|\partial_tU^n\|^2+a_{dG}(U^{n+1},U^{n+1})
\right)
+\delta(\partial_tU^n,U^{n+1}),
\end{equation}
there holds
\begin{equation}\label{4.4}
\frac12\min(1,\nu)E^n(U)
\leq\widetilde{E}^n_\delta(U)
\leq\frac32\max(1,\mu)E^n(U),
\end{equation}
provided 
\begin{equation*}
0<\delta\leq
\min\left\{
\frac{(\alpha\gamma+\beta\nu)}{2\gamma},
\frac{\nu}{2(\beta\nu+\alpha\gamma)}
\right\}.
\end{equation*}
\end{lemma}
\begin{proof}
Choose $\chi_h=\partial_tU^n$ in \eqref{fully1} and 
use
$\partial_t\bar{\partial}_tU^n
=\bar{\partial}_t\partial_tU^n$
and the symmetry of $a_{dG}(\cdot,\cdot)$ to obtain
\begin{eqnarray*}
(\bar{\partial}_t\partial_tU^n,\partial_tU^n)
&=&
\frac12\bar{\partial}_t\|\partial_tU^n\|^2
+\frac{k}{2}\|\bar{\partial}_t\partial_tU^n\|^2,\\
a_{dG}(U^{n+1},\partial_tU^n)
&=&
\frac12\bar{\partial}_t
a_{dG}(U^{n+1},U^{n+1})
+\frac{k}{2}a_{dG}(\partial_tU^n,\partial_tU^n).
\end{eqnarray*}
Hence,
\begin{eqnarray}\label{estimate-d}
\frac12\bar{\partial}_t
\left(
\|\partial_tU^n\|^2+a_{dG}(U^{n+1},U^{n+1})
\right)
&+&\alpha\|\partial_tU^n\|^2
+\beta a_{dG}(\partial_tU^n,\partial_tU^n)
\nonumber\\
&+&\frac{k}{2}\left(
\|\bar{\partial}_t\partial_tU^n\|^2
+a_{dG}(\partial_tU^n,\partial_tU^n)
\right)
=0.
\end{eqnarray}
By \eqref{coereq} and \eqref{L2},
\begin{eqnarray*}
\alpha\|\partial_tU^n\|^2
+\beta a_{dG}(\partial_tU^n,\partial_tU^n)
&\geq&
\vartheta\left(\alpha+\frac{\beta\nu}{\gamma}\right)
\|\partial_tU^n\|^2\\
&&{}+(1-\vartheta)\alpha\|\partial_tU^n\|^2
+(1-\vartheta)\beta\nu\vertiii{\partial_tU^n}^2.
\end{eqnarray*}
Apply the coercivity to the last term on the right hand side of \eqref{estimate-d}. Altogether, we derive \eqref{4.2}.

To verify \eqref{4.4}, we combine inequality \eqref{L2} with the Hölder’s inequality and the Young’s inequality  to 
arrive at
\begin{eqnarray}\label{4.5}
\bigl|(\partial_t U^n, U^{n+1})\bigr|
&\le &\|\partial_t U^{n}\|\,\|U^{n+1}\|
\le \sqrt{\gamma}\|\partial_t U^n\|\,\vertiii{U^{n+1}} \nonumber\\
&\le & \frac{\gamma}{2(\alpha\gamma+\beta\nu)}\|\partial_t U^n\|^2
+\frac{(\alpha\gamma+\beta\nu)}{2}
\vertiii{U^{n+1}}^2.
\end{eqnarray}
Substituting \eqref{4.5} into the definition \eqref{4.3} shows
\begin{eqnarray*}
\widetilde{E}^n_\delta(U)
&\ge &\frac12 (\|\partial_t U^n\|^2 + a_{dG}(U^{n+1},\,U^{n+1}) )
+ \delta(\partial_t U^n, U^{n+1})\nonumber\\
&\ge &\frac14 (\|\partial_t U^n\|^2 + \nu\vertiii{U^{n+1}}^2)
+ \left(\frac{1}{4}
- \frac{\gamma}{2(\alpha\gamma+\beta\nu)}{\delta}\right)\|\partial_t U^n\|^2
+ \left(\frac{\nu}{4}
- \frac{(\alpha\gamma+\beta\nu)}{2}\delta\right)\vertiii{U^{n+1}}^2\nonumber\\
&\ge &\frac{1}{2}E^n(U)\min(1,\,\nu),
\end{eqnarray*}
provided
\[
0 < \delta \le 
\min\!\left (\dfrac{(\alpha \gamma+\beta \nu)}{2\gamma},\,\dfrac{\nu}{2(\beta \nu+\alpha \gamma)}\right ).\]
Similarly, using \eqref{4.5} in \eqref{4.3}, we obtain
\begin{eqnarray*}
\widetilde{E}^n_\delta(U)
&\le &\frac12 (\|\partial_t U^n\|^2 + a_{dG}(U^{n+1},\,U^{n+1}) )
+ \delta(\partial_t U^n, U^{n+1})\nonumber\\
&\le &\frac{3}{4} (\|\partial_t U^n\|^2 +\mu \vertiii{U^{n+1}}^2 )
- \left(\frac{1}{4}
- \frac{\gamma}{2(\alpha\gamma+\beta\nu)}{\delta}\right)\|\partial_t U^n\|^2
- \left(\frac{\mu}{4}
- \frac{(\alpha\gamma+\beta\nu)}{2}\delta\right)\vertiii{U^{n+1}}^2\nonumber\\
&\le &\frac{3}{2}E^n(U)\max(1,\mu), ~\,\text{whenever}\,~ 0 < \delta \le 
\min\!\left (\dfrac{(\alpha \gamma+\beta \nu)}{2\gamma},\,\dfrac{\mu}{2(\beta \nu+\alpha \gamma)}\right ). 
\end{eqnarray*}
 This completes the proof.
\end{proof}
\noindent
The next theorem establishes the uniform exponential decay of the completely discrete energy.
\noindent
\begin{theorem}\label{expdecayfullyap}
For $0 < \delta \le 
\min\!\left (\dfrac{\alpha \gamma+\beta \nu}{2\gamma},\,\dfrac{\nu}{2(\beta \nu+\alpha \gamma)}\right )$, the solution $U^{n+1}$ for $n\ge1$ of
\eqref{fully1} satisfies
\[
E^n(U) \le 3pe^{-\frac{1}{16}\delta t_n}E^0(U),
\]
where
\(
p=\frac{\max(1,\mu)}{\min(1,\nu)}.
\)
\end{theorem}
\begin{proof}
From the definition of the extended energy \eqref{4.3}, there holds
\begin{eqnarray}\label{fullyexp1}
\bar\partial_t \tilde E_\delta^n(U)
&=&\frac12\bar\partial_t \left(
\|\partial_tU^n\|^2+a_{dG}(U^{n+1},U^{n+1})
\right)+\delta \bar\partial_t (\partial_t U^n,U^{n+1}) \nonumber\\
&=&\frac12\bar\partial_t \left(
\|\partial_tU^n\|^2+a_{dG}(U^{n+1},U^{n+1})
\right)
+\delta\|\partial_t U^n\|^2
+\delta(\bar\partial_t \partial_t U^n,U^n).
\end{eqnarray}
Setting $\chi_h=U^n$ in \eqref{fully1} with the relation
$U^n-U^{n+1}=-k\partial_t U^n$, a use of the H\"older inequality with
\eqref{L2}, \eqref{coereq}, the symmetry of $a_{dG}(\cdot,\cdot)$ and the Young inequality imply
\begin{eqnarray}\label{fullyexp2}
(\bar\partial_t \partial_t U^n,U^n)
&=&-\beta a_{dG}(\partial_t U^n,U^n)
-\alpha(\partial_t U^n,U^n)
-a_{dG}(U^{n+1},U^n)\nonumber\\
&=&-\beta a_{dG}(\partial_tU^n,U^{n+1})
+\beta k a_{dG}(\partial_tU^n,\partial_tU^n)-\alpha(\partial_tU^n,U^{n+1})
+\alpha k\|\partial_tU^n\|^2\nonumber\\
&&-a_{dG}(U^{n+1},U^{n+1})
+k a_{dG}(U^{n+1},\partial_tU^n)\nonumber\\
&\leq& \frac{\beta^2}{2\delta_0}
a_{dG}(\partial_tU^n,\partial_tU^n)
+\frac{\delta_0}{2}a_{dG}(U^{n+1},U^{n+1})+\left(\frac{\beta^2}{2\delta_1}
+\frac{k^2\delta_1}{2}\right)
a_{dG}(\partial_tU^n,\partial_tU^n)\nonumber\\ 
&&+\frac{\alpha^2\gamma}{2\nu\delta_0}
\|\partial_tU^n\|^2
+\frac{\delta_0}{2}a_{dG}(U^{n+1},U^{n+1})+\left(
\frac{\alpha^2\gamma}{2\nu\delta_1}
+\frac{k^2\nu\delta_1}{2\gamma}
\right)\|\partial_tU^n\|^2\nonumber\\
&&-a_{dG}(U^{n+1},U^{n+1})+\frac{\delta_2}{2}a_{dG}(U^{n+1},U^{n+1})
+\frac{k^2}{2\delta_2}
a_{dG}(\partial_tU^n,\partial_tU^n)\nonumber\\
&\leq&-\left(1-\delta_0-\frac{\delta_2}{2}\right)
a_{dG}(U^{n+1},U^{n+1})+\frac{\beta^2}{2}
\left(\frac1{\delta_0}+\frac1{\delta_1}\right)
a_{dG}(\partial_tU^n,\partial_tU^n)\nonumber\\
&&+\frac{\alpha^2\gamma}{2\nu}
\left(\frac1{\delta_0}+\frac1{\delta_1}\right)
\|\partial_tU^n\|^2+k^2\left(\delta_1+\frac1{2\delta_2}\right)
a_{dG}(\partial_tU^n,\partial_tU^n).
\end{eqnarray}
Substituting \eqref{fullyexp2} in \eqref{fullyexp1},
using \eqref{4.2}, and choosing
$\vartheta=21/40$, $\delta_0=3/4$, $\delta_1=30/17$,
and $\delta_2=2/5$, we obtain, for
$0<\delta\leq\min\left\{\frac{(\alpha\gamma+\beta\nu)}{(2\gamma)},
\frac{\nu}{2(\alpha\gamma+\beta\nu)}\right\}$
and $\delta k\leq\frac{34}{205}$,
\begin{eqnarray*}
 \bar\partial_t \tilde E_\delta^n(U)&\leq & -\frac{\delta}{20}a_{dG}(U^{n+1},U^{n+1})
-\left(
\frac{21(\alpha\gamma+\beta\nu)}{40\gamma}-\delta
\right)\|\partial_tU^n\|^2\\
&&-\beta\left(
\frac{19}{40}-\frac{19\beta\delta}{20}
\right)a_{dG}(\partial_tU^n,\partial_tU^n)
-\alpha\left(
\frac{19}{40}-\frac{19\alpha\gamma\delta}{20\nu}
\right)\|\partial_tU^n\|^2\\
&&-\frac{k}{2}\|\bar{\partial}_t\partial_tU^n\|^2
-\frac{k}{2}
\left(1-\frac{205}{34}\delta k\right)
a_{dG}(\partial_tU^n,\partial_tU^n)\\
&\leq & -\frac{\delta}{20}
\left(
\|\partial_tU^n\|^2+a_{dG}(U^{n+1},U^{n+1})
\right).
   \end{eqnarray*}
An application of \eqref{L2}, \eqref{coereq} and the Young inequality yields
\begin{eqnarray*}
\delta|(\partial_tU^n,U^{n+1})|
&\leq &
\frac{\delta(\alpha\gamma+\beta\nu)}{2\nu}
\|\partial_tU^n\|^2
+\frac{\delta\gamma}{2(\alpha\gamma+\beta\nu)}
a_{dG}(U^{n+1},U^{n+1})\\
&\leq&
\frac14\left(
\|\partial_tU^n\|^2+a_{dG}(U^{n+1},U^{n+1})
\right).
\end{eqnarray*}
Consequently,
\begin{eqnarray}\label{tildeinequality}
\frac14\left(
\|\partial_tU^n\|^2+a_{dG}(U^{n+1},U^{n+1})
\right)
\leq\widetilde E_\delta^n(U)
\leq
\frac34\left(
\|\partial_tU^n\|^2+a_{dG}(U^{n+1},U^{n+1})
\right).
\end{eqnarray}
Using the upper bound, we conclude that
\begin{equation*}
\bar{\partial}_t\widetilde E_\delta^n(U)
\leq-\frac{\delta}{15}\widetilde E_\delta^n(U).
\end{equation*}
Multiplying the above equation by
 $e^{\upsilon\delta t_n}$ and rearranging gives
\begin{equation}\label{fullyex3}
\frac{e^{\upsilon\delta k}}{k}
\left(
e^{\upsilon\delta t_n}\widetilde E_\delta^n(U)
-e^{\upsilon\delta t_{n-1}}\widetilde E_\delta^{n-1}(U)
\right)+e^{\upsilon\delta t_n}
\left(
\frac{\delta}{15}
-\frac{e^{\upsilon\delta k}-1}{k}
\right)\widetilde E_\delta^n(U)
\leq0.
\end{equation}
Choose $\upsilon=1/16$. Since $\delta k\leq34/205<1$,
the inequality $e^x\leq(1-x)^{-1}$ for $0\leq x<1$
implies
\begin{equation*}
e^{\delta k/16}-1
\leq\frac{\delta k}{16-\delta k}
\leq\frac{\delta k}{15}.
\end{equation*}
Thus the second term  of \eqref{fullyex3} is nonnegative. Dividing by
$e^{\upsilon\delta k}$, multiplying by $k$, and replacing the time index by $j$ and summing from
$j=1$ to $n$, we obtain
\begin{equation*}
\widetilde E_\delta^n(U)
\leq e^{-\upsilon\delta t_n}\widetilde E_\delta^0(U).
\end{equation*}
Utilize \eqref{4.4} to arrive at
\begin{equation}\label{fullyex4}
 E^n(U)
\leq\widetilde{E}^n_\delta(U)\leq e^{-\upsilon\delta t_n}\widetilde E_\delta^0(U)
\leq3pe^{-\upsilon\delta t_n}E^0(U), \quad \text{where} \quad p=\frac{\max(1,\mu)}{\min(1,\nu)}.  \end{equation}
This concludes the rest of the proof.
\end{proof}
\subsection{Fully Discrete Error Estimates}
This subsection is devoted to establish fully discrete error estimates.\\
We split the error at $t_n$ as
\begin{equation*}
e^n:=u(t_n)-U^n
=\bigl(u(t_n)-R_hu(t_n)\bigr)
 +\bigl(R_hu(t_n)-U^n\bigr)
=\xi^n+\rho^n.
\end{equation*}
Subtracting \eqref{fully1} from \eqref{1.5} at time $t_{n+1}$, adding
the appropriate difference quotients and applying orthogonality property of $R_h$  gives
\begin{equation}\label{fullyerreq}
(\bar\partial_t\partial_t\rho^n,\psi_h)
+\beta a_{dG}(\partial_t\rho^n,\psi_h)
+\alpha(\partial_t\rho^n,\psi_h)
+a_{dG}(\rho^{n+1},\psi_h)
=(\Upsilon^n,\psi_h)
+\beta a_{dG}(\sigma_5^n,\psi_h),
\end{equation}
for every $\psi_h\in V_h$ and $n\geq1$, where 
$
\Upsilon^n:=\sum_{i=1}^{4}\sigma_i^n,\,
\sigma_1^n:=-\bar\partial_t\partial_t\xi^n,\,
\sigma_2^n:=-\alpha\partial_t\xi^n,
\sigma_3^n:=\bar\partial_t\partial_tu^n-u''(t_{n+1}),\,
\sigma_4^n:=-\alpha\bigl(u'(t_{n+1})-\partial_tu^n\bigr),\,
\sigma_5^n:=-\bigl(u'(t_{n+1})-\partial_tu^n\bigr).
$
For convenience, we define
\begin{equation*}
\mathcal{D}_r
:=
\beta^2\|u_1\|_{r+3}^2
+\|u_1\|_{r+2}^2
+\|u_0\|_{r+3}^2.
\end{equation*}
For a sufficiently smooth function $w$, the forward
and central differences satisfy
\begin{eqnarray}
\partial_tw^n
&=&\frac1k\int_{t_n}^{t_{n+1}}w'(s)\,ds,\label{forward}\\
\bar\partial_t\partial_tw^n
&=&\frac1{k^2}\int_{-k}^{k}
(k-|s|)w''(t_n+s)\,ds,\qquad \frac1{k^2}\int_{-k}^{k}(k-|s|)\,ds=1.\label{central}
\end{eqnarray}
A use of the Cauchy--Schwarz inequality,  \eqref{approx1}, \eqref{semiw13} and \eqref{central} imply
\begin{eqnarray}\label{sigma1n}
\|\sigma_1^n\|^2
&\leq&
\frac1{k^2}\int_{-k}^{k}
(k-|s|)\|\xi''(t_n+s)\|^2\,ds\\
&\leq&
\frac{Ch^{2r+2}}{k^2}\int_{-k}^{k}
(k-|s|)\|u''(t_n+s)\|_{r+1}^2\,\leq Ch^{2r+2}e^{-2\delta t_{n-1}}\mathcal{D}_r.
\end{eqnarray}
Similarly, from \eqref{approx1} and \eqref{semiw16}, \eqref{forward}, we have 
\begin{eqnarray}\label{sigma2n}
\|\sigma_2^n\|^2
&\leq&
\frac{\alpha^2}{k}
\int_{t_n}^{t_{n+1}}\|\xi'(s)\|^2\,ds
\leq
\frac{C\alpha^2h^{2r+2}}{k}
\int_{t_n}^{t_{n+1}}\|u'(s)\|_{r+1}^2\,ds\\&\leq& Ch^{2r+2}e^{-2\delta t_n}
\mathcal{D}_r.
\end{eqnarray}
The integral representation of the central difference \eqref{central}
gives
\begin{equation*}
\sigma_3^n
=\frac1{k^2}\int_{-k}^{k}(k-|s|)
\bigl(u''(t_n+s)-u''(t_{n+1})\bigr)\,ds.
\end{equation*}
Since
\begin{equation*}
u''(t_n+s)-u''(t_{n+1})
=-\int_{t_n+s}^{t_{n+1}}u'''(q)\,dq,
\end{equation*}
we have
\begin{eqnarray*}
\|\sigma_3^n\|
&\leq&
\frac1{k^2}\int_{-k}^{k}
(k-|s|)(k-s)\,ds\,
\max_{t_{n-1}\leq q\leq t_{n+1}}\|u'''(q)\|\\
&=&
k\max_{t_{n-1}\leq q\leq t_{n+1}}\|u'''(q)\|.
\end{eqnarray*}
To control the third order time derivative, differentiate
\eqref{1.1} with respect to time
\begin{equation*}
u'''=\beta \Delta u''-\alpha u''+\Delta u'.
\end{equation*}
Hence,
\begin{equation*}
\|u'''(t)\|^2
\leq C\left(
\beta^2\|u''(t)\|_2^2
+\alpha^2\|u''(t)\|^2
+\|u'(t)\|_2^2
\right).
\end{equation*}
Since $r+1\geq2$, the regularity estimates
\eqref{semiw15} and \eqref{semiw16} yield
\begin{equation*}
\|u'''(t)\|^2
\leq Ce^{-2\delta t}\mathcal D_r.
\end{equation*}
Consequently,
\begin{equation}\label{sigma3n}
\|\sigma_3^n\|^2
\leq Ck^2e^{-2\delta t_{n-1}}\mathcal D_r.
\end{equation}
An integration in time gives
\begin{eqnarray*}
\sigma_5^n
=\frac1k\int_{t_n}^{t_{n+1}}
\bigl(u'(s)-u'(t_{n+1})\bigr)\,ds
=-\frac1k\int_{t_n}^{t_{n+1}}
(s-t_n)u''(s)\,ds.
\end{eqnarray*}
Thus, for $m=0,1$,
\begin{equation*}
\|\sigma_5^n\|_m
\leq\frac{k}{2}
\max_{t_n\leq s\leq t_{n+1}}\|u''(s)\|_m.
\end{equation*}
In particular, since $\sigma_4^n=\alpha\sigma_5^n$, we apply \eqref{semiw15} to obtain
\begin{eqnarray}
\|\sigma_4^n\|^2
&\leq&
\frac{\alpha^2k^2}{4}
\max_{t_n\leq s\leq t_{n+1}}\|u''(s)\|^2
\leq C\alpha^2k^2e^{-2\delta t_n}\mathcal D_r,
\label{sigma4n}\\
|\sigma_5^n|_1^2
&\leq&
\frac{k^2}{4}
\max_{t_n\leq s\leq t_{n+1}}|u''(s)|_1^2
\leq Ck^2e^{-2\delta t_n}\mathcal D_r.\label{sigma5n}
\end{eqnarray}
Choosing $\chi_h=\partial_t\rho^n$ in
\eqref{fullyerreq} and using the symmetry of $a_{dG}(\cdot,\cdot)$, \eqref{L2} and similar techniques used in Lemma \ref{fullyD1}
yield
\begin{eqnarray}\label{fullyerrtest}
\frac12\bar\partial_t
\left(
\|\partial_t\rho^n\|^2
+a_{dG}(\rho^{n+1},\rho^{n+1})
\right)
&+&
\vartheta\left(\frac{\alpha\gamma+\beta\nu}{\gamma}\right)
\|\partial_t\rho^n\|^2
+(1-\vartheta)\alpha\|\partial_t\rho^n\|^2
\nonumber\\&+&(1-\vartheta)\beta
a_{dG}(\partial_t\rho^n,\partial_t\rho^n)
+\frac{k}{2}
\left(
\|\bar\partial_t\partial_t\rho^n\|^2
+a_{dG}(\partial_t\rho^n,\partial_t\rho^n)
\right)\nonumber\\
&\le &Ck^{-1}(\|\Upsilon^n\|^2+|\sigma_5^n|_1)+\frac{k}{4}a_{dG}(\partial_t\rho^n,\,\partial_t\rho^n).
\end{eqnarray}
From the definition of the extended energy \eqref{4.3} and using the arguments applied in Theorem \ref{expdecayfullyap}, we have
\begin{equation}\label{fullyerrtest1}
\bar\partial_t \tilde E_\delta^n(\rho)
=\frac12\bar\partial_t \left(
\|\partial_t\rho^n\|^2+a_{dG}(\rho^{n+1},\rho^{n+1})
\right)
+\delta\|\partial_t \rho^n\|^2
+\delta(\bar\partial_t \partial_t \rho^n,\rho^n).
\end{equation}
Substituting $\psi_h= \rho^n$ in \eqref{fullyerreq} with the relation
$\rho^n-\rho^{n+1}=-k\partial_t \rho^n$, a use of the H\"older inequality with
\eqref{L2}, \eqref{coereq}, the symmetry of $a_{dG}(\cdot,\cdot)$ and the Young inequality yield
\begin{eqnarray}\label{fullyerrtest2}
(\bar\partial_t\partial_t\rho^n, \rho^n)&=&
-\beta a_{dG}(\partial_t\rho^n, \rho^n)
-\alpha(\partial_t\rho^n, \rho^n)
-a_{dG}(\rho^{n+1}, \rho^n)
+(\Upsilon^n, \rho^n)
+\beta a_{dG}(\sigma_5^n, \rho^n)\nonumber\\
&\leq& -\left(1-\delta_0-\frac{\delta_2}{2}\right)
a_{dG}(U^{n+1},U^{n+1})+\frac{\beta^2}{2}
\left(\frac1{\delta_0}+\frac1{\delta_1}\right)
a_{dG}(\partial_t\rho^n,\partial_t\rho^n)\nonumber\\
&&+\frac{\alpha^2\gamma}{2\nu}
\left(\frac1{\delta_0}+\frac1{\delta_1}\right)
\|\partial_t\rho^n\|^2+k^2\left(\delta_1+\frac1{2\delta_2}\right)
a_{dG}(\partial_t\rho^n,\partial_t\rho^n)\nonumber\\
&&+C\left(\|\Upsilon^n\|^2+|\sigma_5^n|_1^2\right)+ \frac{\delta_2}2a_{dG}(\rho^n,\rho^n)+\frac{k^2}{2\delta_2}
a_{dG}(\partial_tU^n,\partial_tU^n).
\end{eqnarray}
Using \eqref{fullyerrtest}  and \eqref{fullyerrtest2}  in \eqref{fullyerrtest1} with \eqref{tildeinequality},
 and choosing
$\vartheta=\frac{21}{40}$, $\delta_0=\frac{3}{4}$, $\delta_1=\frac{30}{17}$,
and $\delta_2=\frac{1}{5}$, we obtain the following estimate, provided that $0<\delta\leq\min\left\{\frac{(\alpha\gamma+\beta\nu)}{(2\gamma)},
\frac{\nu}{2(\alpha\gamma+\beta\nu)}\right\}$
and $\delta k\leq\frac{17}{460}$:
\begin{eqnarray*}
  \bar\partial_t \tilde E_\delta^n(\rho)&\leq &-\frac{1}{10} \delta  \left(\|\partial_t\rho^n\|^2
+a_{dG}(\rho^{n+1},\rho^{n+1})\right)-\frac{k}{2}\left(\frac{1}{2}-\frac{230}{17}\delta k\right)a_{dG}(\partial_tU^n,\partial_tU^n)\\
& &+Ck^{-1}(\|\Upsilon^n\|^2+|\sigma_5^n|_1^2)\\
&\leq& -\frac{1}{15} \delta\tilde E_\delta^n(\rho) +Ck^{-1}(\|\Upsilon^n\|^2+|\sigma_5^n|_1^2).
\end{eqnarray*}
Since $\delta k\leq17/460<1$, applying the analogous techniques to the above equation that utilized for \eqref{fullyex4}, we find that
\begin{equation}\label{fullyerrweightedbound}
E^n(\rho)
\leq
3p\,e^{-\delta t_n/16}E^0(\rho)
+C\sum_{j=1}^{n}
e^{-\delta(t_n-t_{j-1})/16}
\left(
|\Upsilon^j|^2+|\sigma_5^j|_1^2
\right).
\end{equation}
 For the derivation of $E^0(\rho)$,  compare \eqref{fullyerrstartup} with the exact
equation at $t_1=k$ and using the projection
orthogonality gives
\begin{eqnarray}\label{fullyerrinitialeq}
\frac{2}{k}(\partial_t\rho^0,\psi_h)
+\beta a_{dG}(\partial_t\rho^0,\psi_h)
&+&\alpha(\partial_t\rho^0,\psi_h)
\nonumber\\
&+&a_{dG}(\rho^1,\psi_h)
=(\Upsilon^0,\psi_h)
+\beta a_{dG}(\sigma_5^0,\psi_h),\quad \psi_h \in V_h.
\end{eqnarray}
where $\Upsilon^0=\sum_{i=1}^{4}\sigma_i^0$ and $
\sigma_1^0:=-\frac{2}{k}\partial_t\xi^0,
\quad
\sigma_2^0:=-\alpha\partial_t\xi^0,\quad
\sigma_3^0:=
\frac{2}{k}(\partial_tu(0)-u_1)-u''(t_1),\quad
\sigma_4^0:=
\alpha\bigl(\partial_tu(0)-u'(t_1)\bigr),\quad
\sigma_5^0:=
\partial_tu(0)-u'(t_1).
$
Twice applications of the fundamental theorem of calculus, followed
by a change in the order of integration, give
\begin{eqnarray*}
u(t_1)-u(0)-ku_1
=\int_0^{t_1}(t_1-q)u''(q)\,dq.
\end{eqnarray*}
It follows that
\begin{eqnarray*}
\sigma_3^0
=\frac{2}{k^2}\int_0^{t_1}(t_1-s)
\bigl(u''(s)-u''(t_1)\bigr)\,ds
=-\frac{2}{k^2}\int_0^{t_1}(t_1-s)
\int_s^{t_1}u'''(q)\,dq\,ds.
\end{eqnarray*}
Taking norms and using the triangle inequality yields
\begin{eqnarray*}
\|\sigma_3^0\|
\leq\frac{2}{k^2}\int_0^{t_1}(t_1-s)^2\,ds
\max_{0\leq q\leq t_1}\|u'''(q)\|
=\frac{2k}{3}
\max_{0\leq q\leq t_1}\|u'''(q)\|.
\end{eqnarray*}
Applying similar techniques used in $\|\Upsilon^n\|$ and $\sigma_5^n$ for the remaining terms of $\Upsilon^0$ and $\sigma_5^0$, we have
\begin{equation}\label{fullyerrinitialforcing}
\|\Upsilon^0\|+|\sigma_5^0|_1
\leq C\left(\frac{h^{r+1}}{k}+k\right)
\mathcal{D}_r^{1/2}.
\end{equation}
Since $\rho^0=0$, we have $\rho^1=k\partial_t\rho^0$.
Choosing $\psi_h=k\partial_t\rho^0$ in
\eqref{fullyerrinitialeq} gives
\begin{eqnarray*}
2\|\partial_t\rho^0\|^2
+a_{dG}(\rho^1,\rho^1)
+\beta k
a_{dG}(\partial_t\rho^0,\partial_t\rho^0)
+\alpha k\|\partial_t\rho^0\|^2
=k(\Upsilon^0,\partial_t\rho^0)
+\beta a_{dG}(\sigma_5^0,\rho^1).
\end{eqnarray*}
Applying \eqref{coereq}, Cauchy--Schwarz, Young's inequalities and \eqref{fullyerrinitialforcing}, we obtain
\begin{eqnarray}\label{fullyerrinitialforcing1}
2\min\{2,\nu\}E^0(\rho^0)\leq 2\|\partial_t\rho^0\|^2
+\nu \vertiii{\rho^1}^2
\leq k(\Upsilon^0,\partial_t\rho^0)
+\beta a_{dG}(\sigma_5^0,\rho^1)\leq C\left(h^4+k^2\right)\mathcal{D}_r.
\end{eqnarray}
By utilizing  \eqref{sigma1n}--\eqref{sigma5n} with \eqref{fullyerrinitialforcing1} in \eqref{fullyerrweightedbound}, we arrive at 
\begin{eqnarray}\label{fullyerrweightedbound10}
E^n(\rho)
&\leq &
C\left(h^4+k^2\right)\,e^{-\delta t_n/16}\mathcal{D}_r
+C\left(h^4+k^2\right)\sum_{j=1}^{n}
e^{-\delta(t_n-t_{j-1})/16}e^{-2\delta t_{j-1}}\nonumber\\
& \leq &C\left(h^4+k^2\right)\,e^{-\delta t_n/16}\mathcal{D}_r.
\end{eqnarray}
The main theorem now follows from the error decomposition,
the triangle inequality, and the estimates in \eqref{approx2},  \eqref{fullyerrweightedbound10} and Theorem~\ref{regularity}.
\begin{theorem}\label{fullyerrtheorem}
Let $U^0=R_hu_0$, let $U^1$ be determined by
\eqref{fullyerrstartup}, and let $U^{n+1}$ satisfy
\eqref{fully1} for $n\geq1$ and 
\(0<\delta\leq
\min\left\{
\frac{\alpha\gamma+\beta\nu}{2\gamma},
\frac{\nu}{2(\alpha\gamma+\beta\nu)}
\right\},
\)
Then, We have 
\begin{eqnarray*}
\vertiii{u(t_n)-U^n}^2
&\leq&
C\left(
h^{r}+k^2\right)e^{-\delta t_n/16}\mathcal D_r,\\
\|u(t_n)-U^n\|^2
&\leq&
C\left(
h^{r+1}+k^2\right)e^{-\delta t_n/16}\mathcal D_r,
\end{eqnarray*}
where
\begin{equation*}
\mathcal{D}_r
:=
\beta^2\|u_1\|_{r+3}^2
+\|u_1\|_{r+2}^2
+\|u_0\|_{r+3}^2.
\end{equation*}
The constant $C$ is independent of $h$, $k$, and $n$.
\end{theorem}
\section{Some Generalizations}\label{s5}
In this section, we discuss  extension of our analysis to nonhomogeneous strongly damped wave equations and to a class of strongly  damped fourth order  wave equations. However, we only indicate the difference in the analysis for both the problems.

\subsection{ Nonhomogeneous Problems}
For the strongly damped wave equation with non zero forcing function $f$:
\begin{equation*}
u''-\beta\Delta u'+\alpha u'-\Delta u=f
\quad\text{in }\Omega\times(0,\infty),
\end{equation*}
with boundary and initial conditions \eqref{1.2}--\eqref{1.3},
where $f\in L^2(\Omega)$ is independent of time.
Let $u_\infty\in H_0^1(\Omega)$ be the unique weak
solution of
\begin{equation*}
-\Delta u_\infty=f\quad\text{in }\Omega,
\qquad
u_\infty=0\quad\text{on }\partial\Omega.
\end{equation*}
We now study the convergence to this equilibrium by introducing
\(
v(t):=u(t)-u_\infty.
\) 
Note that $v$ satisfies 
\begin{equation}\label{inhomogeneous-v}
v''-\beta\Delta v'+\alpha v'-\Delta v=0
\quad\text{in }\Omega\times(0,\infty),
\end{equation}
with boundary \eqref{1.2} by replacing $u$ by $v$  and initial conditions:
\begin{equation*}
 v(x,0) = u_0(x)-u_{\infty}(x)\;\;\;\mbox{and}\;\;\; v'(x,0)=u_1(x),\;\;x\in \Omega.  
\end{equation*}
Note  the results of Theorem ~ \ref{regularity} are valid also for the solution $v$ of \eqref{inhomogeneous-v}.
With corresponding Galerkin approximation $u_h(t) \in V_h$ for each $t>0$ satisfying 
\begin{equation}\label{semiwin}
(u_h'',\phi_h)+\beta a_{dG}(u_h',\phi_h)
+\alpha(u_h',\phi_h)+a_{dG}(u_h,\phi_h)
=(f,\phi_h),
\qquad \phi_h\in V_h,
\end{equation}
we  define the discrete equilibrium $u_{\infty,h}\in V_h$
by
\begin{equation}\label{discretesteadystate}
a_{dG}(u_{\infty,h},v_h)=(f,v_h),
\qquad v_h\in V_h.
\end{equation}
A use of coercivity  and boundedness shows using Lax-Milgram theorem that there exists a unique solution$u_{\infty,h}\in V_h$ of  \eqref{coereq}. 
Setting
\(
v_h(t):=u_h(t)-u_{\infty,h},
\)
and subtracting \eqref{discretesteadystate} from \eqref{semiwin} yields
\begin{equation*}
(v_h'',\phi_h)+\beta a_{dG}(v_h',\phi_h)
+\alpha(v_h',\phi_h)+a_{dG}(v_h,\phi_h)=0,
\qquad \phi_h\in V_h.
\end{equation*}
Consequently, the semidiscrete decay
estimates apply to $v_h$, with initial data
\(
v_h(0)=u_{0,h}-u_{\infty,h},
\quad
v_h'(0)=u_{1,h}.
\)
 As in Theorem ~\ref{semidecay} for $u_h$, similar estimate is valid for $u_h$ replaced by $v_h.$ In particular, under the  decay conditions
\(
0<\delta\leq
\frac{\alpha\gamma+\beta\nu}{2\gamma},
\quad
\delta<\frac{\nu}{\alpha\gamma+\beta\nu},
\)
we arrive at
\begin{equation*}
E_h(v_h)(t)
\leq Ce^{-2\delta t}E_h(v_h)(0).
\end{equation*}
Now apply Theorems~\ref{semierr2}--\ref{semierr3} of  Section~\ref{s3} by replacing $u_h$ by $v_h$ to arrive at
\begin{equation*}
    h\vertiii{(v-v_h)(t)} + \|(v-v_h)(t)\| = O \big(e^{-\delta t}\;h^{r+1}\big).
\end{equation*}

\begin{remark}
Now as in Theorem~\ref{semimax}, we easily derive 
that  the semidiscrete solution converges to its
equilibrium at the exponential rate $\delta$:
\begin{equation*}
\|u_h(t)-u_{\infty,h}\|_{L^\infty(\Omega)}
=O(e^{-\delta t})
\qquad\text{as }t\to\infty.
\end{equation*}
\end{remark}
\begin{remark}
Suppose that $f(t)=O(e^{-\delta_0t})$ and that
the homogeneous solution decays at rate $\delta$
in the corresponding energy norm. If
$\delta\neq\delta_0$, the forced solution decays
at rate $\delta_*=\min(\delta,\delta_0)$.
When $\delta=\delta_0$, the general estimate is
$O(e^{-\delta t})$, which implies exponential
decay at every rate $0<\delta_*<\delta$.
\end{remark}

\subsection{ Fourth Order Evolution Equations}
This section discusses the extension of our analysis to a class of fourth-order strongly damped wave equations. The initial-boundary value problem is to find a  function $u : \Omega \times [0, T] \to \mathbb{R}$ such that
\begin{equation} \label{bihar1}
    u_{tt}  + \beta \Delta^2 u_t + \alpha u_t + \Delta^2 u =0 \quad \text{in } \Omega \times (0, T],
\end{equation}
subject to the homogeneous clamped boundary conditions:
\begin{equation} \label{bihar2}
    u(x,t) = 0, \quad \frac{\partial u}{\partial n}(x,t) = 0 \quad \text{on } \partial\Omega \times [0, T],
\end{equation}
and the prescribed initial conditions:
\begin{equation} \label{bihar3}
    u(x,0) = u_0(x), \quad u'(x,0) = u_1(x) \quad \text{in } \Omega.
\end{equation}
Here, $n$ denotes the outward unit normal to the boundary $\partial \Omega$,  the damping parameters $\alpha >0$ and $\beta>0.$ 

We now introduce the following discontinuous function space:
\begin{equation*}
X := \big\{\phi\in L^2(\Omega):\ \phi|_T\in W^{4,\,4/3}(T)\ \text{for all } T\in\mathcal{T}_h \big\}.
\end{equation*}
Define the symmetric bilinear form
$a_{dG}^b:X\times X\to\mathbb R$ by
\begin{equation}\label{bihar5}
a_{dG}^b(\phi,v)
:=a_{\mathrm{pw}}(\phi,v)+b_h(\phi,v)+c_{dG}(\phi,v),
\qquad \phi,v\in X.
\end{equation}
The same notation is used for its restriction to
$X_h\times X_h$, since $X_h\subset X$. The discontinuous Galerkin (dG) energy norm $\vertiii{\cdot}_{\text{dG}}$ on the space $X$ is defined by
\begin{equation*} 
\vertiii{{v}}_{\text{dG}} := \left( \vertiii{v}_{\text{pw}}^2 + c_{\text{dG}}(v, v) \right)^{1/2} \quad \forall v \in X,
\end{equation*}
where $\vertiii{v}_{\text{pw}}^2 = a_{\text{pw}}(v, v)$ denotes the piecewise volumetric energy.
\\
We shall frequently employ the following $L^2$-estimate (see \cite[Lemma 2.4]{georgoulis2013}): For any
$\phi\in X$
\begin{equation}\label{birL2}
\|\phi\|^2_{L^2(\Omega)} \le \gamma_b\,\vertiii{\phi}_{dG}^2,
\end{equation}
where  $\gamma_b>0$ is a constant.
Assume that the exact solution satisfies
$u(t),u'(t)\in H_0^2(\Omega)\cap W^{4,4/3}(\Omega)$
and $u''(t)\in L^2(\Omega)$.
Then elementwise integration by parts yields
\begin{equation}\label{bi1.5}
(u'',v)+\beta a_{dG}^b(u',v)
+\alpha(u',v)+a_{dG}^b(u,v)=0,
\qquad \forall v\in X,\quad t>0.
\end{equation}

In the dG FEM method, we define the fully discontinuous piecewise polynomial  finite element space $V_h$ as
\begin{equation*} 
    X_h =\mathbb{P}_r(\mathcal{T}_h),\, r\ge 2.
\end{equation*}

Let $\mathcal{E}_h$ denote the set of all edges in the triangulation $\mathcal{T}_h$. We introduce the discrete bilinear form $a_{dG}^b: V_h \times V_h \to \mathbb{R}$, which is decomposed as
\begin{equation} 
    a_{dG}^b(\phi_h, v_h) := a_{\text{pw}}(\phi_h, v_h) + b_h(\phi_h, v_h) + c_{\text{dG}}(\phi_h, v_h).
\end{equation}
The components of the above equation  are defined as follows. The volumetric Hessian term is given by
\begin{equation*} 
    a_{\text{pw}}(\phi_h, v_h) := \sum_{K \in \mathcal{T}_h} \int_{K} D_{pw}^2 \phi_h : D_{pw}^2 v_h \,\mathrm{d}x,
\end{equation*}
where $D_{pw}^2$ denotes the piecewise Hessian matrix and $:$ represents the Frobenius inner product of matrices. The consistency and symmetry term is defined on the mesh skeleton by
\begin{eqnarray*}
    b_h(\phi_h, v_h) &=& \sum_{e \in \mathcal{E}_h} \int_e [\phi_h] \{\mathbf{n}_e \cdot \nabla \Delta v_h\} \, \mathrm{d}s 
    + \sum_{e \in \mathcal{E}_h} \int_e [v_h] \{\mathbf{n}_e \cdot \nabla \Delta \phi_h\} \, \mathrm{d}s \nonumber \\
    && - \sum_{e \in \mathcal{E}_h} \int_e [\nabla \phi_h] \cdot \{D_{\mathrm{pw}}^2 v_h\}\mathbf{n}_e \, \mathrm{d}s 
    - \sum_{e \in \mathcal{E}_h} \int_e [\nabla v_h] \cdot \{D_{\mathrm{pw}}^2 \phi_h\}\mathbf{n}_e \, \mathrm{d}s.
\end{eqnarray*}
Finally, the discontinuous Galerkin interior penalty term is defined as
\begin{equation*} 
    c_{\text{dG}}(\phi_h, v_h) := \sum_{e \in \mathcal{E}_h} \frac{\sigma_{\text{dG}}^1}{h_e^3} \int_e [\phi_h][v_h] \, \mathrm{d}s + \sum_{e \in \mathcal{E}_h} \frac{\sigma_{\text{dG}}^2}{h_e} \int_e \left[ \frac{\partial \phi_h}{\partial n_e} \right] \left[ \frac{\partial v_h}{\partial n_e} \right] \, \mathrm{d}s,
\end{equation*}
where $\sigma_{\text{dG}}^1, \sigma_{\text{dG}}^2 > 0$ are the stabilization penalty parameters. \\
In the $C^0$ interior penalty ($C^0$-IP) method, we restrict the discrete space to globally continuous, piecewise quadratic functions:
\begin{equation*}
    X_h := \mathbb{P}_k(\mathcal{T}_h) \cap H_0^1(\Omega),\, k\ge 2.
\end{equation*}
Because functions in $X_h$ are continuous across element interfaces, the stabilization term only penalizes jumps in the normal derivatives. The discrete bilinear form is defined as
\begin{eqnarray*}
    a_{dG}^{b}(w_h, v_h) &:=& a_{\text{pw}}(w_h, v_h) + b_h^{IP}(w_h, v_h) + c_{\text{IP}}(w_h, v_h),\nonumber \\ 
a_{dG}(\phi_h,\psi_h) &:=&
\sum_{T\in\mathcal{T}_h}\int_T  \nabla\phi_h\cdot \nabla\psi_h\,d\bx, \nonumber \\ 
b_h^{IP}(\phi_h, v_h) &:=& - \sum_{e \in \mathcal{E}_h} \int_e [\nabla \phi_h] \cdot \{D_{\text{pw}}^2 v_h\}n_e \, \mathrm{d}s - \sum_{e \in \mathcal{E}_h} \int_e [\nabla v_h] \cdot \{D_{\text{pw}}^2 \phi_h\}n_e \, \mathrm{d}s,
\end{eqnarray*}
where $a_{\text{pw}}$ and $b_h$ are the standard volumetric and consistency terms, and the boundary penalty is
\begin{equation*}
    c_{\text{IP}}(w_h, v_h) := \sum_{e \in \mathcal{E}_h} \frac{\sigma_{\text{IP}}}{h_e} \int_e \left[ \frac{\partial w_h}{\partial n_e} \right] \left[ \frac{\partial v_h}{\partial n_e} \right] \, \mathrm{d}s,
\end{equation*}
with $\sigma_{\text{IP}} > 0$. The associated energy norm on $V_h$ is defined by
\begin{equation*}
    \vertiii{v_h}_{\text{IP}} := \left( \vertiii{v_h}_{\text{pw}}^2 + c_{\text{IP}}(v_h, v_h) \right)^{1/2}.
\end{equation*}
The semi-discrete DG formulation seeks a discrete solution $u_h(t) \in X_h$ such that
\begin{equation} \label{bihar9}
    (u_{htt}, v_h) + \beta a_{dG}^b(u_{ht}, v_h)+ \alpha (u_{ht}, v_h) + a_{dG}^b(u_h, v_h) =0 \quad \forall v_h \in X_h
\end{equation}
with  initial conditions
\begin{equation} \label{eq:discrete_initial}
    u_h(0) = u_{0h}, \quad u_{ht}(0) = u_{1h},
\end{equation}
where $u_{0h}$ and $u_{1h}$ in $X_h$ are approximations of the exact initial data $u_0$ and $u_1$ to be defined later.

Assume sufficiently large penalty parameters so that
the symmetric bilinear forms \cite[Lemma 5]{Huang2026} satisfy
\begin{eqnarray}\label{biharcoercivity}
a_{dG}^b(v_h,v_h)
&\geq&\nu_b\vertiii{v_h}_{dG}^2,\\
\label{biharcontinuity}
|a_{dG}^b(v_h,w_h)|
&\leq&\mu_b\vertiii{v_h}_{dG}\vertiii{w_h}_{dG},
\end{eqnarray}
for all $v_h,w_h\in X_h$, where the positive constants
are independent of $h$.
Define the discrete energy by
\begin{equation*}
E_h^b(u_h)(t)
:=\frac12\left(
\|u_h'(t)\|^2
+\vertiii{u_h(t)}_{dG}^2
\right).
\end{equation*}
\begin{theorem}\label{bihardecay}
Let $u_h$ satisfy \eqref{bihar9} and
initial conditions \eqref{eq:discrete_initial}.
Further, let
\begin{equation*}
0<\delta\leq
\min\left\{
\frac{\alpha\gamma_b+\beta\nu_b}{2\gamma_b},
\frac{\nu_b}{\beta\nu_b+\alpha\gamma_b}
\right\}.
\end{equation*}
Then the following estimates hold for all $t\geq0$.
\begin{itemize}
\item[(i)] If
\(
\delta<\frac{\nu_b}{\beta\nu_b+\alpha\gamma_b},
\)
then
\begin{equation*}
E_h^b(u_h)(t)
\leq Ce^{-2\delta t}E_h^b(u_h)(0).
\end{equation*}
This includes
$\delta=(\alpha\gamma_b+\beta\nu_b)/(2\gamma_b)$
whenever this value is strictly below
$\nu_b/(\beta\nu_b+\alpha\gamma_b)$.

\item[(ii)] If
\(
\delta=\frac{\nu_b}{\beta\nu_b+\alpha\gamma_b},
\)
then
\begin{equation*}
\|u_h'(t)\|^2
\leq Ce^{-2\delta t}E_h^b(u_h)(0).
\end{equation*}
\end{itemize}
The constant $C$ is independent of $h$ and $t$.
\end{theorem}

\begin{proof}
Choose $u_{h,\delta}=e^{\delta t}u_h$  in   \eqref{bihar9} and apply \eqref{delta_derivative} to obtain
\begin{eqnarray*}
(u_{h,\delta}'',v_h)
+\beta a_{dG}^b(u_{h,\delta}',v_h)
&+&(\alpha-2\delta)(u_{h,\delta}',v_h)
+(1-\beta\delta)a_{dG}^b(u_{h,\delta},v_h)\nonumber\\
&+&
-\delta(\alpha-\delta)(u_{h,\delta},v_h)=0,
\qquad v_h\in X_h.
\end{eqnarray*}
With $v_h=u_{\delta}'$, rewrite it as
\begin{equation*} 
\frac{d}{dt} \mathcal {F}_{h,\delta}^b(t) +  \beta a_{dG}^b(u_{h,\delta}',u_{h,\delta}')+ (\alpha-2\delta)(u_{h,\delta}',\,u_{h,\delta}')= 0,
\end{equation*}
where
\begin{eqnarray*}
\mathcal F_{h,\delta}^b(t)
&:=&\frac12\Big(
\|u_{h,\delta}'\|^2
+(1-\beta\delta)
a_{dG}^b(u_{h,\delta},u_{h,\delta})
-\delta(\alpha-\delta)\|u_{h,\delta}\|^2
\Big).
\end{eqnarray*}
Now, introduce 
\begin{eqnarray*}
\widetilde E_{h,\delta}^b(t)
&:=&\frac12\Big(
\|u_{h,\delta}'\|^2
+(1-\beta\delta)\nu_b\vertiii{u_{h,\delta}}^2
-\delta(\alpha-\delta)\|u_{h,\delta}\|^2
\Big).
\end{eqnarray*}
Replacing $\mathcal F_{h,\delta}^b(t)$ and
$\widetilde E_{h,\delta}^b(t)$ by
$\mathcal F_{h,\delta}(t)$ and
$\widetilde E_{h,\delta}(t)$, respectively, and following
the argument in the proof of Theorem~\ref{semidecay},
together with \eqref{biharcoercivity} and
\eqref{biharcontinuity}, completes the proof.
\end{proof}
We recall the following polynomial approximation result from \cite{brenner2008mathematical}.
\begin{lemma}
\label{bihar-local-approx}
Let $w\in H^s(K)$, where $s\ge2$, and set
$\ell=\min\{s,r+1\}$. There exists a polynomial
$I_Kw\in\mathbb P_r(K)$ such that
\begin{equation*}
|w-I_Kw|_{j,K}
\le C h_K^{\ell-j}|w|_{\ell,K},
\qquad 0\leq j \leq \ell.
\end{equation*}
The constant $C$ is independent of $h_K$ and $w.$
\end{lemma}
Now, define the piecewise polynomial operator $I_h$ by
\(
(I_hw)|_K:=I_K(w|_K),
\qquad K\in\mathcal T_h.
\)
\begin{lemma}\cite{Huang2026}
\label{bihar-global-approx}
Let $w\in H^s(\Omega)\cap H_0^2(\Omega)$,
where $s\ge2$, and let $\ell=\min\{s,r+1\}$.
Then,
\begin{equation}\label{bihar-global-estimate}
\vertiii{w-I_hw}_{dG}
\le Ch^{\ell-2}|w|_{\ell,\Omega}.
\end{equation}
The constant $C$ is independent of $h$ and $w$.
\end{lemma}

Below, we indicate the error estimates which are uniformly exponentially decaying. 
Let $M_h u$ be the projection satisfying
\begin{equation}\label{mh:projection}
 a_{dG}^b(u-M_hu,v_h)=0
\qquad v_h\in X_h,
\end{equation}
Differentiating \eqref{mh:projection}, there holds
\begin{eqnarray}\label{mh:projection1}
 a_{dG}^b((u-M_hu)', v_h)=0,\quad M_h u(0)= I_h u_0
\qquad v_h\in X_h, \\
a_{dG}^b((u-M_hu)'', v_h)=0,\quad M_h u'(0)= I_h u_1
\qquad v_h\in X_h.\label{mh:projection2}
\end{eqnarray}
Let $e=u-M_hu+M_hu-u_h=\xi_b+\rho_b$. Using \eqref{mh:projection}--\eqref{mh:projection2}, Lemmas~\ref{bihar-local-approx} and~\ref{bihar-global-approx} and 
and adapting the arguments in Theorems~3 and~4
of \cite{Huang2026}, we obtain the following results.  Set
\(
q_r:=r-1+\min\{r-1,2\}, \quad r\ge2.
\)
For the  projection $M_h$, the standard approximation estimates  yield
\begin{eqnarray}
\sum_{m=0}^{2}
\left(
\|\xi_b^{(m)}(t)\|
+h^{q_r-r+1}\vertiii{\xi_b^{(m)}(t)}_{dG}
\right)
&\leq&
Ch^{q_r}
\sum_{m=0}^{2}\|u^{(m)}\|_{r+3}.
\label{biapprox1}\\
\sum_{m=0}^{2}
\left(
\|\xi_{b,\delta}^{(m)}(t)\|
+h^{q_r-r+1}\vertiii{\xi_{b,\delta}^{(m)}(t)}_{dG}
\right)
&\leq&
Ch^{q_r}e^{\delta t}
\sum_{m=0}^{2}\|u^{(m)}\|_{r+3}.\nonumber
\end{eqnarray}
The constant $C>0$ is independent of $h$ and $t$.
The error equation follows by subtracting \eqref{bihar9}
from \eqref{bi1.5} and applying the orthogonality property
of the fourth order elliptic projection:
\begin{equation}\label{mh:weightedprojection}
(\rho_b'', \phi_h) + \beta\, a_{dG}^b(\rho_b', \phi_h) + \alpha(\rho_b', \phi_h)
+ a_{dG}^b(\rho_b, \phi_h)
= -(\xi''+\alpha\xi', \phi_h),
\quad \forall\; \phi_h \in X_h .
\end{equation}
Setting $\rho_{b,\delta}=e^{\delta t}\rho_b$ and multiplying
\eqref{mh:weightedprojection} by $e^{\delta t},$  we arrive at
\begin{eqnarray}\label{mh:erroreq}
(\rho_{b,\delta}'', \phi_h) &+& \beta\, a_{dG}(\rho_{b,\delta}', \phi_h)
+ (\alpha - 2\delta)(\rho_{b,\delta}', \phi_h)
+ (1 - \beta\delta)\, a_{dG}(\rho_{b,\delta}, \phi_h)\nonumber\\
& -& \delta(\alpha - \delta)(\rho_{b,\delta}, \phi_h)
= -e^{\delta t}(\xi'', \phi_h) - \alpha e^{\delta t}(\xi', \phi_h), \quad \forall\;  \phi_h \in X_h .
\end{eqnarray}
Now, define
\begin{equation*}
E_h^b(\rho_b)(t):=\frac12\left(
\|\rho_b'(t)\|^2
+\vertiii{\rho_b(t)}_{dG}^2\right).
\end{equation*}
\begin{theorem}\label{mh:theorem33}
Let $\rho$ satisfy \eqref{mh:erroreq}, and let
$0<\delta_0<\alpha\gamma_b+\beta\nu_b$. Suppose
\begin{eqnarray}\label{mh:deltarange}
0<\delta\le\min\left\{
\frac{\alpha\gamma_b+\beta\nu_b-\delta_0}{2\gamma_b},
\frac{\nu}{\beta\nu_b+\alpha\gamma_b}\right\}.
\end{eqnarray}
If $\delta<\nu_b/(\beta\nu_b+\alpha\gamma_b)$, then, for every $t\ge 0$,
\begin{eqnarray*}
E_h^b(\rho_b)(t)
\le Ce^{-2\delta t}E_h^b(\rho_b)(0)+C\int_0^t e^{-2\delta(t-s)}
\left(\|\xi_b''(s)\|^2+\alpha^2\|\xi_b'(s)\|^2\right)\,ds.
\end{eqnarray*}
Here, $\delta=\nu/(\beta\nu+\alpha\gamma)$, provided
\eqref{mh:deltarange} also holds and 
the constant $C$ is independent of $h$ and $t$.
\end{theorem}
\begin{proof}
Choose $\rho_{b,\delta}=e^{\delta t}\rho_b$. Multiply 
\eqref{mh:erroreq} by $e^{\delta t}$  to find that
\begin{align*}
(\rho_{b,\delta}'',v_h)&+\beta a_{dG}^b(\rho_{b,\delta}',v_h)
+(\alpha-2\delta)(\rho_{b,\delta}',v_h)\nonumber\\
&+(1-\beta\delta)a_{dG}^b(\rho_{b,\delta},v_h)
-\delta(\alpha-\delta)(\rho_{b,\delta},v_h)
=-e^{\delta t}(\xi_b''+\alpha\xi_b',v_h).
\end{align*}
With $v_h=\rho_{b,\delta}'$, rewrite it as
\[
\frac{d}{dt} F_{h,\delta}^{b}(t)
+\beta a_{dG}^b(\rho_{b,\delta}',\rho_{b,\delta}')
+(\alpha-2\delta)\|\rho_{b,\delta}'\|^2
=-e^{\delta t}(\xi_b''+\alpha\xi_b',\rho_{b,\delta}'),
\]
where
\begin{equation*}
 F_{h,\delta}^{b}(t):=\frac12\Big(
\|\rho_{b,\delta}'\|^2
+(1-\beta\delta)a_{dG}^b(\rho_{b,\delta},\rho_{b,\delta})
-\delta(\alpha-\delta)\|\rho_{b,\delta}\|^2\Big).
\end{equation*}
The proof follows by replacing $F_{h,\delta}^b(t)$ with
$\mathcal F_{h,\delta}(t)$ and adapting the argument of
Theorem~\ref{semidecay}, using \eqref{biharcoercivity}
and \eqref{biharcontinuity}. This concludes the rest of the proof.
\end{proof}
Let $u_{0h}=M_hu_0$, and choose $u_{1h}$ as the
$L^2$-projection or a suitable interpolant of $u_1$.
Then $\rho_b(0)=0$ and $\rho_b'(0)=M_hu_1-u_{1h}$;
hence,
\[
E_h^b(\rho_b)(0)\le Ch^{2q_r}\|u_1\|_{r+3}^2,
\qquad q_r=r-1+\min\{r-1,2\}.
\]
Under the required regularity, compatibility, and decay
assumptions, the continuous decay estimates with
$j=(r+5)/4$ and $k=1,2$, together with
$u''(0)=-\beta\Delta^2u_1-\alpha u_1-\Delta^2u_0$,
yield
\[
\sum_{m=0}^{2}\|u^{(m)}(t)\|_{r+3}^2
\le Ce^{-2\delta t}\mathcal D_b,
\qquad
\mathcal D_b:=
\beta^2\|u_1\|_{r+7}^2
+\|u_1\|_{r+5}^2+\|u_0\|_{r+7}^2.
\]
Combining these bounds with Theorem~\ref{mh:theorem33}
and \eqref{biapprox1}, and using
$\int_0^t e^{-2\delta(t-s)}e^{-2\delta s}\,ds
=te^{-2\delta t}$, gives
\[
\|\rho_b'(t)\|^2+\vertiii{\rho_b(t)}_{dG}^2
\le Ch^{2q_r}(1+t)e^{-2\delta t}\mathcal D_b.
\]
The triangle inequality and the projection estimates
therefore imply
\[
\|u'(t)-u_h'(t)\|^2
+h^{2(q_r-r+1)}\vertiii{u(t)-u_h(t)}_{dG}^2
\le Ch^{2q_r}(1+t)e^{-2\delta t}\mathcal D_b.
\]
For $\beta=0$ and $\alpha>0$, the same argument applies
with $\mathcal D_b=\|u_1\|_{r+5}^2+\|u_0\|_{r+7}^2$
and the corresponding admissible decay rate.
\begin{theorem}\label{bihar:semierr3}
Let $u$ and $u_h$ denote the solutions of
\eqref{bihar1} and \eqref{bihar9}, respectively.
Let $0<\delta_0<\alpha\gamma_b+\beta\nu_b$ and
\(
0<\delta\le
\frac{\alpha\gamma_b+\beta\nu_b-\delta_0}{2\gamma_b},
\qquad
\delta<
\frac{\nu_b}{\beta\nu_b+\alpha\gamma_b}.
\)
Then, for every $t\ge0$,
\begin{equation*}
\|u(t)-u_h(t)\|
\le Ch^{q_r}(1+t)^{1/2}e^{-\delta t}
\left(\|u_0\|_{r+5}+\|u_1\|_{r+3}\right),
\end{equation*}
where $q_r=r-1+\min\{r-1,2\}$.
The constant $C$ is independent of $h$ and $t$.
\end{theorem}
\begin{proof}
Integrating \eqref{mh:weightedprojection} with respect to time, we obtain
\begin{align}\label{bir:semiw20}
(\rho',\phi_h)
+ \beta a_{dG}(\rho,\phi_h)
+ \alpha(\rho,\phi_h)
+ a_{dG}(\hat{\rho},\phi_h)
&=
-(\xi'+\alpha\xi,\phi_h)
\nonumber\\
&\quad
+(\alpha e(0)+e'(0),\phi_h)
+\beta a_{dG}(\rho(0),\phi_h).
\end{align}
We choose $u_{0h}$ and $u_{1h}$ as the elliptic projection and the $L^2$--projection of
$u_0$ and $u_1$, respectively. With this choice, $(e'(0),\phi_h)=0$ and $\rho(0)=0$.\\

\noindent
Let $v=\rho+\tilde{\delta}\hat{\rho}$ in \eqref{bir:semiw20}. Using the definition of the
energy  together with Young’s inequality, \eqref{L2},  \eqref{biharcoercivity}, and \eqref{biharcontinuity}, we arrive at
\begin{equation*}
\frac{d}{dt}\widetilde{E}^b_{\tilde{\delta}}(\rho)(t)
+
F^b(\rho)(t)
\le
C\big(
\|\xi'(t)\|^{2}
+\|\xi(t)\|^{2}
+\|\xi(0)\|^{2}
\big),
\end{equation*}
where
\[
\widetilde{E}^b_{\tilde{\delta}}(\rho)(t)
=
\frac12\Big(
\|\rho\|^{2}
+
a_{dG}^b(\hat{\rho}, \hat{\rho})
\Big)
+
\tilde{\delta}(\rho,\hat{\rho})
+
\frac{\beta\tilde{\delta}}{2}a_{dG}^b(\hat{\rho},\, \hat{\rho})
+
\frac{\alpha\tilde{\delta}}{2}\|\hat{\rho}\|^{2},
\]
with
\[
E^b(\rho)(t)
=
\frac12\Big(
\|\rho(t)\|^{2}
+
\vertiii{\hat{\rho}(t)}_{dG}^{2}
\Big),
\]
and
\[
F^b(\rho)(t)
=
(\tilde{\delta}-\delta_0)a_{dG}^b(\hat{\rho},\,\hat{\rho})
+
(\alpha+\frac{\beta \nu}{\gamma}-\delta_0-\tilde{\delta})\|\rho\|^{2}.
\]
The conclusion follows by adapting the argument used in the proof of
Theorem~\ref{semierr3} together with \eqref{birL2}, \eqref{biharcoercivity}, \eqref{biharcontinuity}.
\end{proof}
\begin{remark}
The preceding analysis extends to the $C^0$ interior
penalty formulation defined above, with
$\vertiii{\cdot}_{dG}$ replaced by
$\vertiii{\cdot}_{IP}$ and $\gamma_b$, $\nu_b$, and
$\mu_b$ denoting the corresponding method-dependent
constants. For the error analysis, $I_h$ is replaced
by a globally continuous interpolant, and $M_h$ is
defined by \eqref{mh:projection} on the $C^0$ interior
penalty space. Under the stated regularity assumptions
and the projection estimate \eqref{biapprox1}, the same
decay and error arguments apply.
\end{remark}
\section{Numerical Experiments}\label{s6}
This section is devoted to numerical experiments that corroborate our theoretical findings. The following examples confirm the predicted convergence rates of the proposed scheme and demonstrate its stability and accuracy in approximating the strongly damped wave equation.

\subsection{Second Order Strongly  Damped wave Equations}
We assess the accuracy and energy decay of the DG approximation
using the homogeneous problem \eqref{1.1}--\eqref{1.3}.
Examples~\ref{ex1}--\ref{example3} use the fully discrete
scheme described in Section~\ref{s4}. The computations are
performed in \texttt{FreeFem++} on successively refined meshes,
with the polynomial degree specified in each example.
All simulations are performed up to $T=1$, with the time step
chosen as $k=O(h^{r+1})$ when using discontinuous piecewise
polynomials of degree $r\geq1$.\\
For errors $E_{h_i}$ and $E_{h_{i+1}}$ measured in the same
norm on two successive meshes, the observed convergence
rate is calculated as
\[
\mathrm{Rate}
=
\frac{\log(E_{h_i}/E_{h_{i+1}})}
     {\log(h_i/h_{i+1})}.
\]

\begin{example}\label{ex1}
Consider the unit square $\Omega=(0,1)^2$ with damping
parameters $\alpha=\pi$ and $\beta=1/\pi$.
The function
\[
u(x,y,t)=e^{-\pi t}\sin(\pi x)\sin(\pi y)
\]
satisfies the homogeneous equation and vanishes on
$\partial\Omega$. Its initial displacement and velocity are
\[
u(x,y,0)=\sin(\pi x)\sin(\pi y),
\qquad
u'(x,y,0)=-\pi\sin(\pi x)\sin(\pi y).
\]
We approximate this solution using discontinuous $P_1$  with penalty parameters $\sigma_e=10$.
\end{example}

\begin{figure}[!ht]
    \centering
    \begin{subfigure}[b]{0.48\textwidth}
        \centering
        \includegraphics[width=\textwidth]{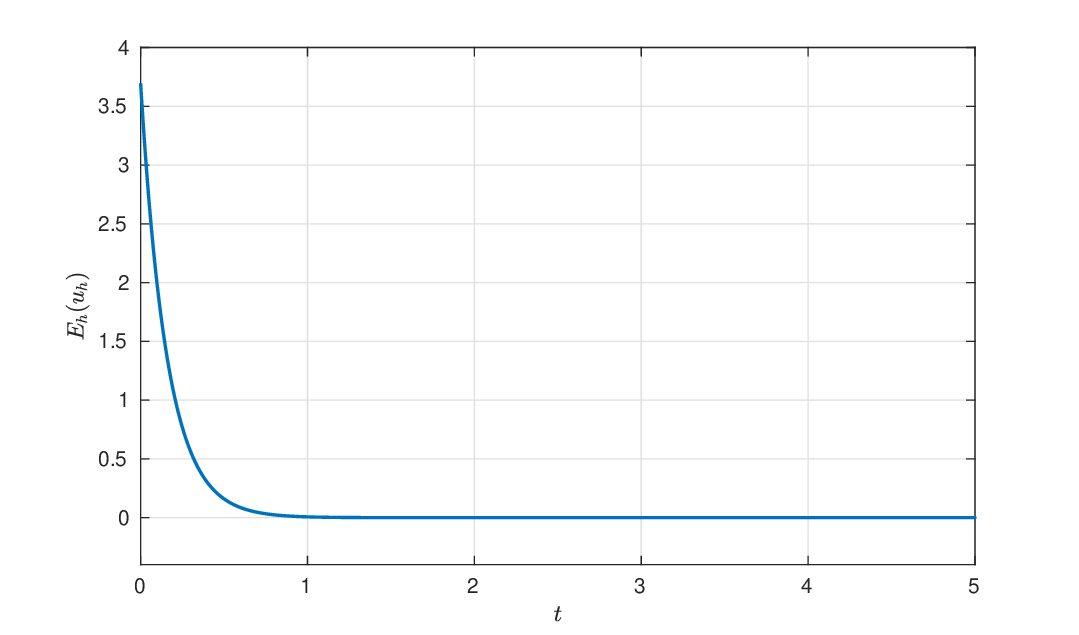}
        \caption{}
        \label{fig:discrete_energyex1}
    \end{subfigure}
    \hfill
    \begin{subfigure}[b]{0.48\textwidth}
        \centering
        \includegraphics[width=\textwidth]{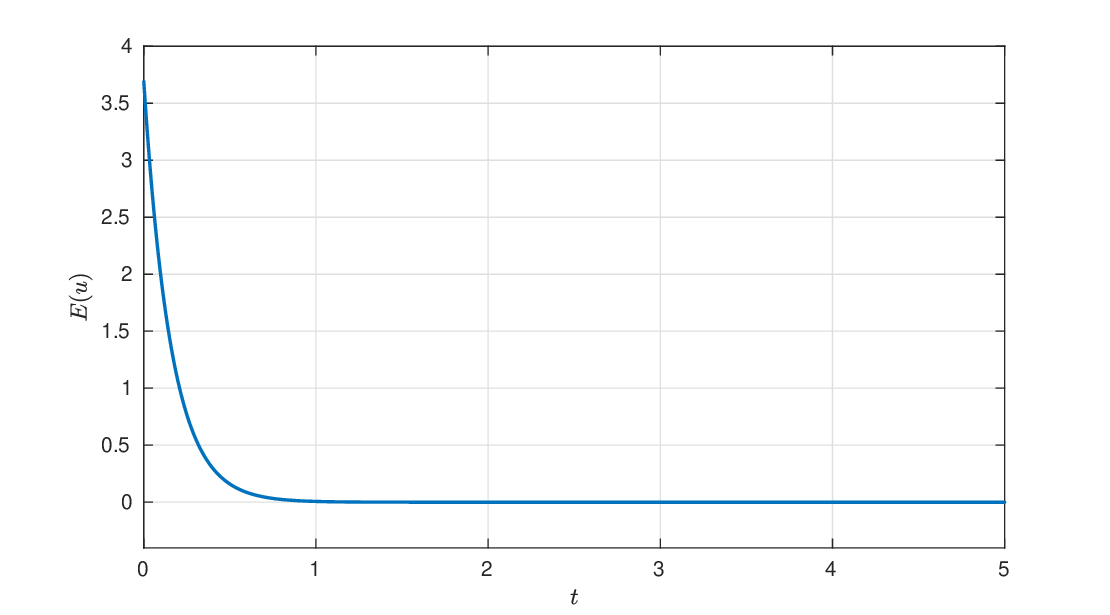}
        \caption{}
        \label{fig:continuous_energyex01}
    \end{subfigure}
\caption{Discrete and continuous energy histories, $E_h(u_h)$ and $E(u)$ for Example~\ref{ex1}.}
    \label{fig:continuous_energyex0}
    \end{figure}

\begin{figure}[!ht]
\center
\includegraphics[scale=0.8]{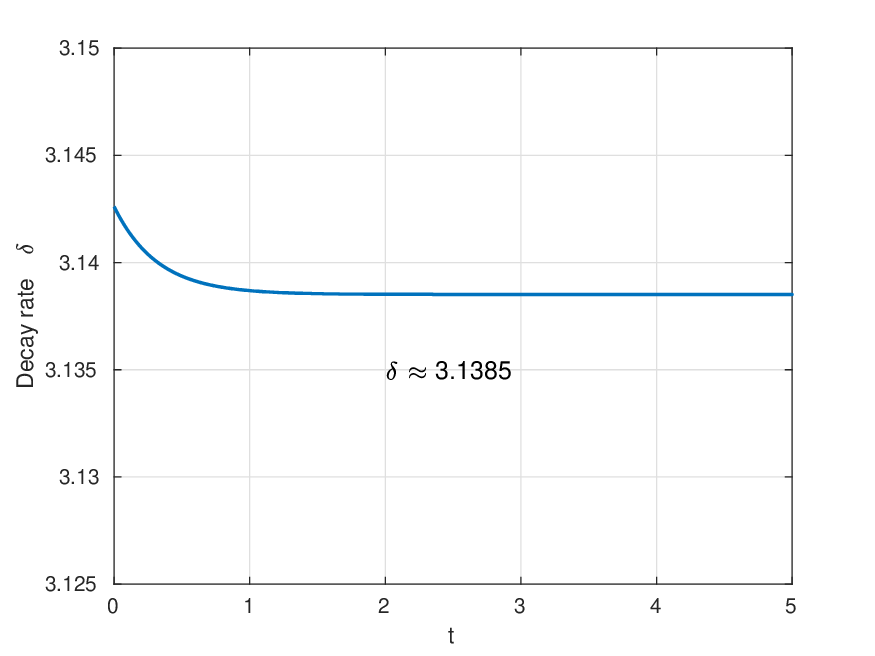}
\caption{ Estimated discrete energy decay rates for Example~\ref{ex1}.}
\label{fig:decayrateex1}
\end{figure}

Table~\ref{errorex1} lists the errors and convergence
rates in the $L^2$, $L^\infty$, and energy norms.
The reported rates agree with the theoretical predictions.
Figure~\ref{fig:continuous_energyex0} compares the discrete
energy $E_h(u_h)$ with the continuous energy $E(u)$ over
the simulation interval, while
Figure~\ref{fig:decayrateex1} displays the computed
exponential decay rates.
\begin{table}[h!]
\centering
\begin{tabular}{|c|cc|cc|cc|}
\hline
\multirow{2}{*}{$N$} 
& \multicolumn{6}{c|}{$\alpha=\pi,\quad \beta=\frac{1}{\pi}$} \\ \cline{2-7}
& $\|u-U\|_{L^2}$ & Rate 
& $\|u-U\|_{L^\infty}$ & Rate 
& $\vertiii{u-U}$ & Rate \\
\hline
10  
& 1.8290(-4) & --    
& 3.5793(-4) & --    
& 4.02227(-2) & -- \\
\hline
20 
& 5.0169(-5) & 1.8662 
& 1.0203(-4) & 1.8107 
& 2.00507(-2) & 1.0043 \\
\hline
30 
& 2.2817(-5) & 1.9432 
& 4.5516(-5) & 1.9908 
& 1.34050(-2) & 0.9930 \\
\hline
40 
& 1.2946(-5) & 1.9698 
& 2.6010(-5) & 1.9451 
& 1.00784(-2) & 0.9914 \\
\hline
50 
& 8.3401(-6) & 1.9724 & 1.6832(-5) & 1.9506 & 8.07763(-3) & 0.9917 \\
\hline
\end{tabular}
\caption{Errors at $T=1$ and convergence rates
for  $P_1$ DG approximation in Example~\ref{ex1}.}
\label{errorex1}
\end{table}

\begin{figure}[!ht]
    \centering
    \begin{subfigure}[b]{0.48\textwidth}
        \centering
        \includegraphics[width=\textwidth]{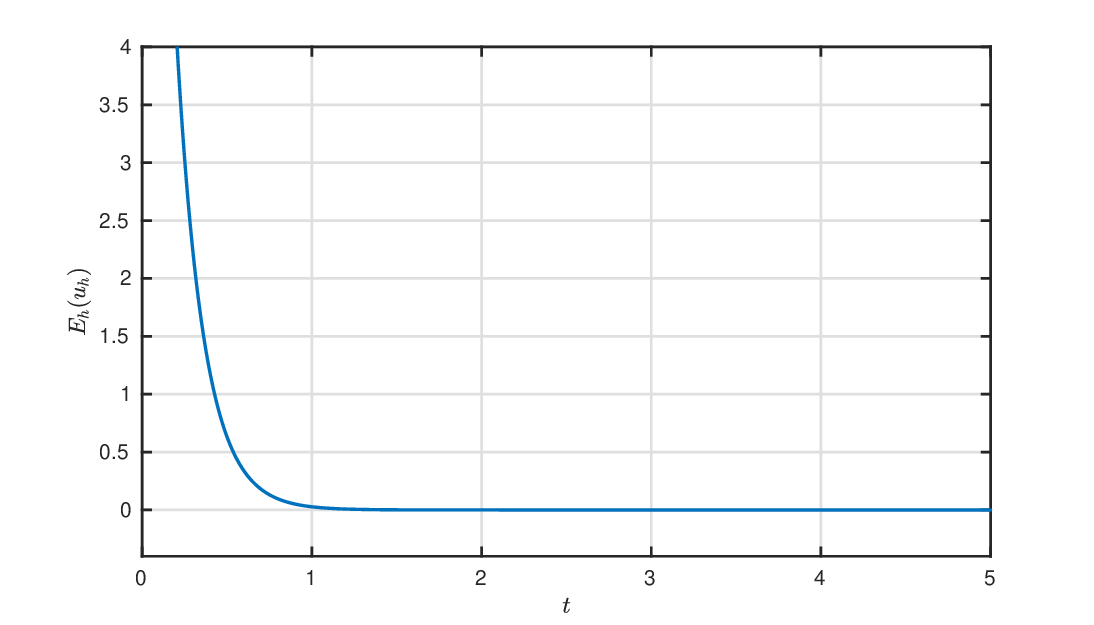}
        \caption{}
        \label{fig:discrete_energyex2}
    \end{subfigure}
    \hfill
    \begin{subfigure}[b]{0.48\textwidth}
        \centering
        \includegraphics[width=\textwidth]{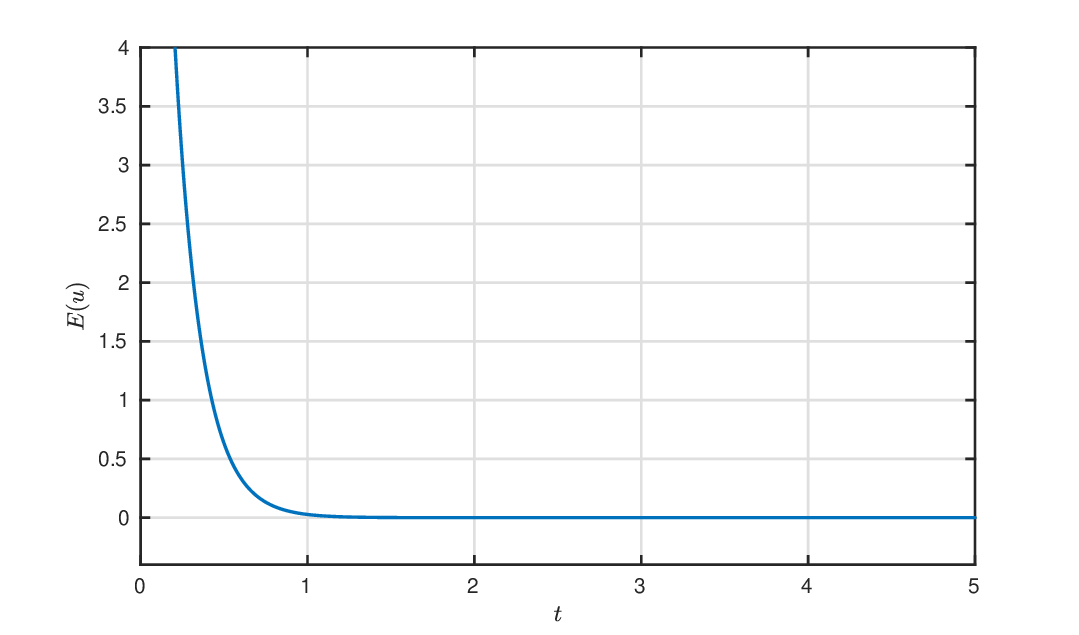}
        \caption{}
        \label{fig:continuous_energyex21}
    \end{subfigure}
    \caption{Comparison of the discrete energy $E_h(u_h)$ and continuous energy $E(u)$, with computed decay rates for the Example~\ref{ex2}.}
    \end{figure}

\begin{figure}[!ht]
\center
\includegraphics[scale=0.8]{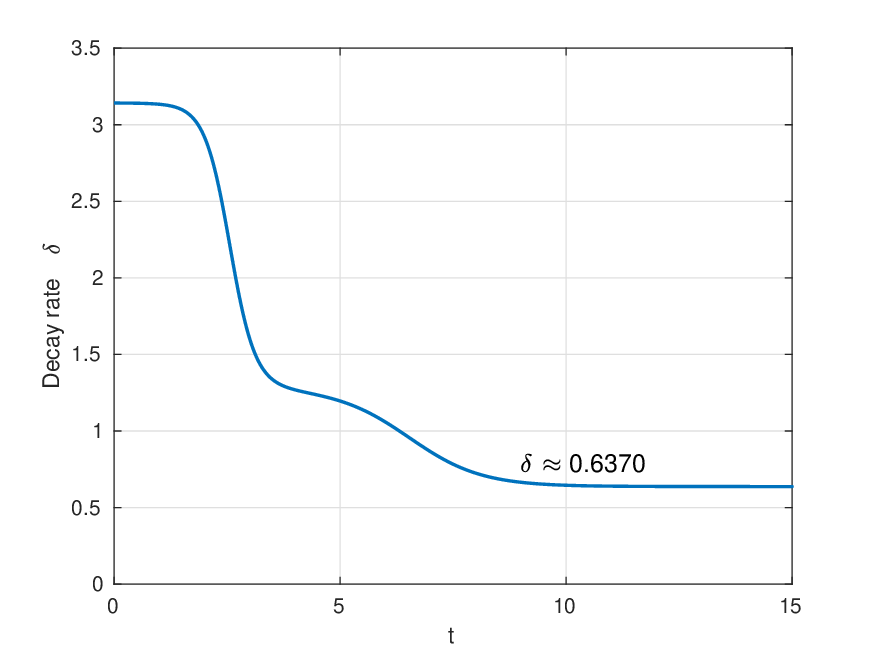}
\caption{ Discrete energy decay proﬁles and estimate of decay rate Example \ref{ex2}.}
\label{fig:decayrate}
\end{figure}

\begin{figure}[!ht]
    \centering
    \begin{subfigure}[b]{0.48\textwidth}
        \centering
        \includegraphics[width=\textwidth]{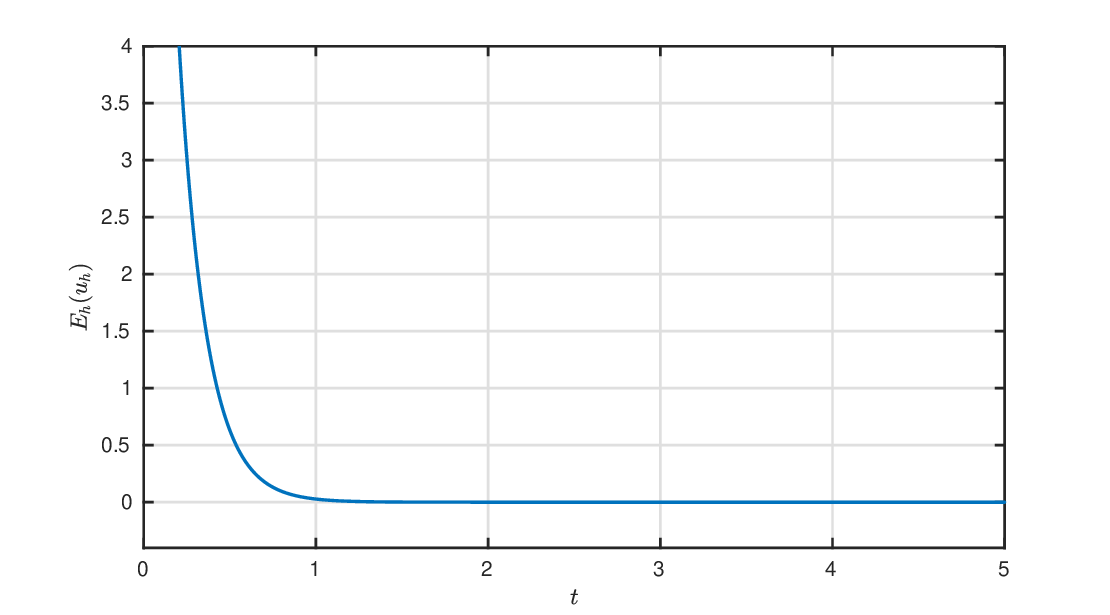}
        \caption{}
        \label{fig:discrete_energy}
    \end{subfigure}
    \hfill
    \begin{subfigure}[b]{0.48\textwidth}
        \centering
        \includegraphics[width=\textwidth]{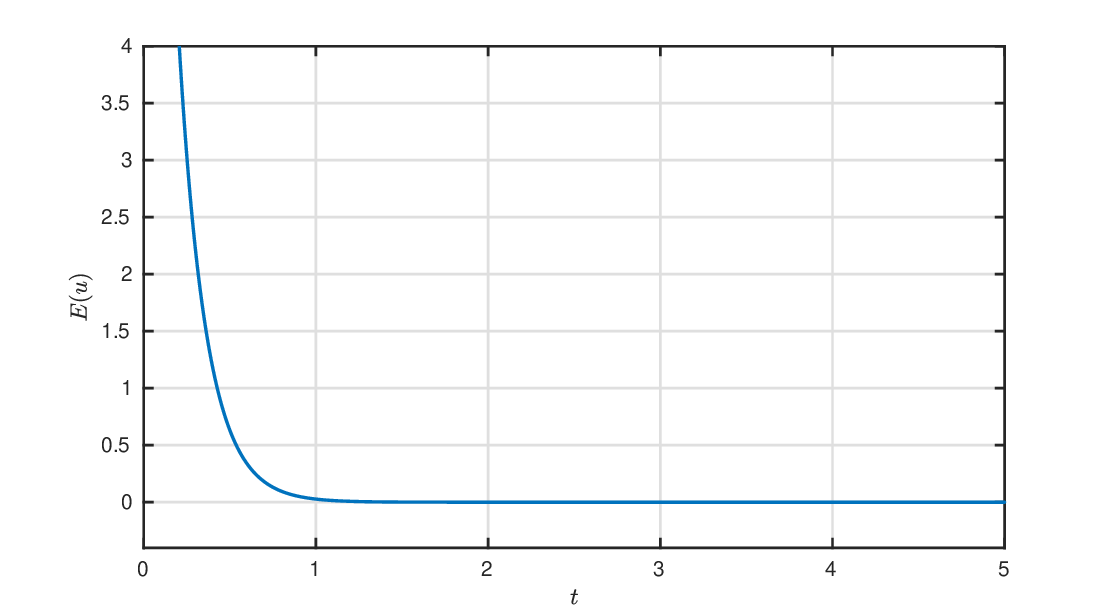}
        \caption{}
        \label{fig:continuous_energy}
    \end{subfigure}
    \caption{Discrete and continuous energy evolution for the strong-damping case in Example~\ref{example3}(i).}
    \label{fig:energy_decay_comparison}
\end{figure}

\begin{figure}[!ht]
\center
\includegraphics[scale=0.8]{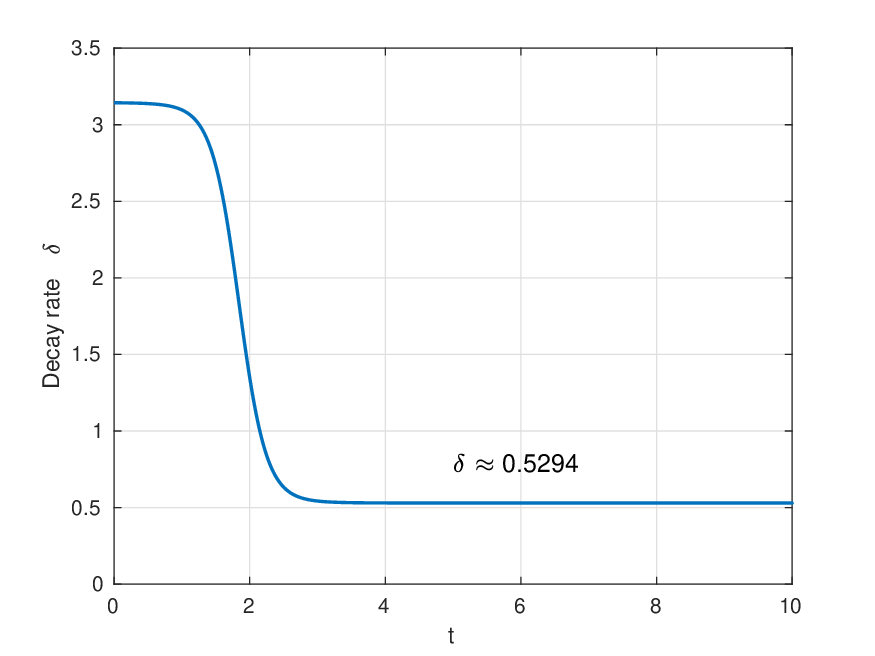}
\caption{Discrete energy  estimated discrete energy decay rates for the Example~\ref{example3}(i).}
\label{fig:decayratealpha0}
\end{figure}

\begin{figure}[!ht]
    \centering
    \begin{subfigure}[b]{0.48\textwidth}
        \centering
        \includegraphics[width=\textwidth]{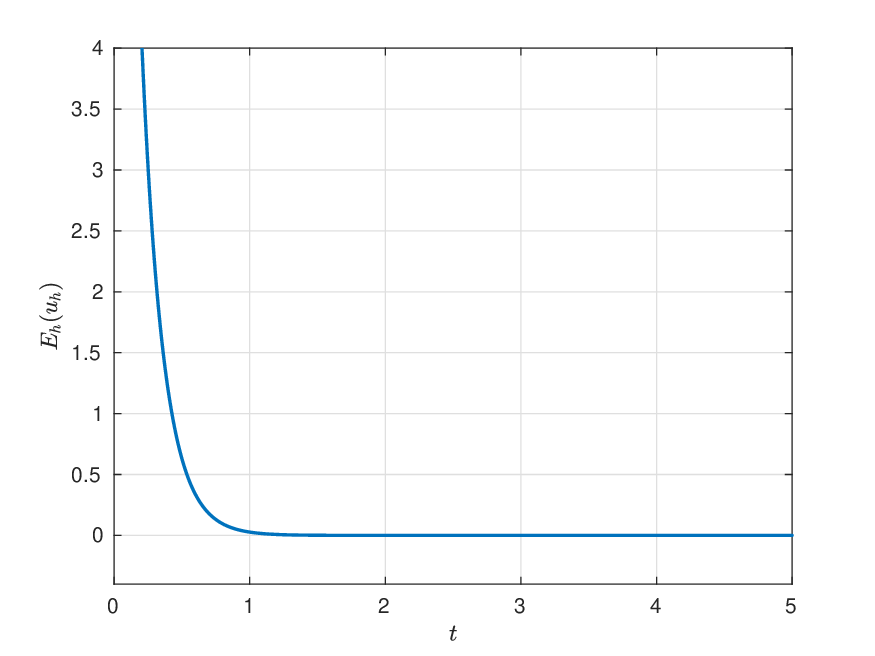}
        \caption{}
        \label{fig:discrete_energy0}
    \end{subfigure}
    \hfill
    \begin{subfigure}[b]{0.48\textwidth}
        \centering
    \includegraphics[width=\textwidth]{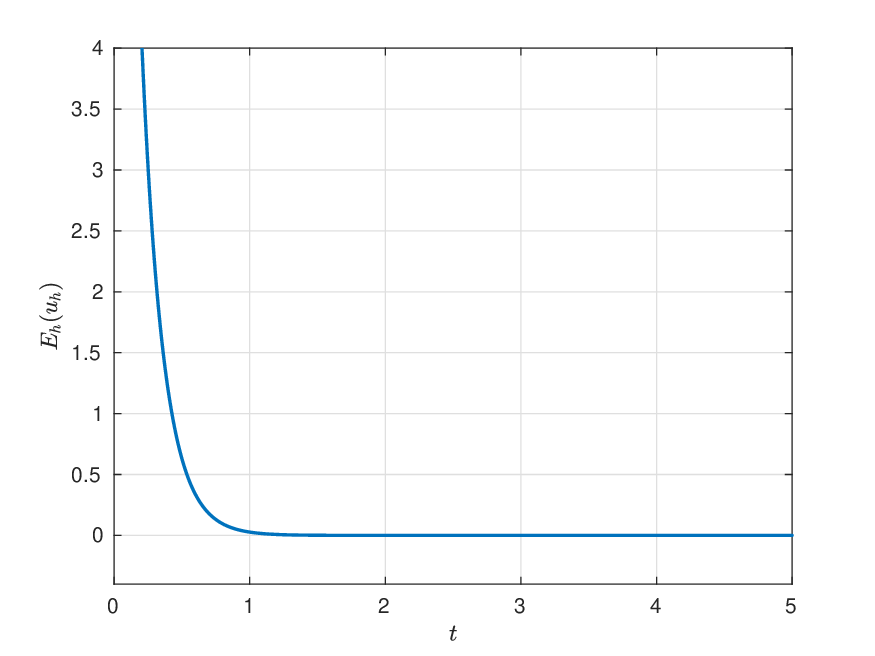}
        \caption{}
        \label{fig:continuous_energy1}
    \end{subfigure}
    \caption{Discrete and continuous energy histories for the weak-damping case in Example~\ref{example3}(ii).}
    \label {fig:energy_decay_comparisonex3}
\end{figure}

\begin{figure}[!ht]
\center
\includegraphics[scale=0.8]{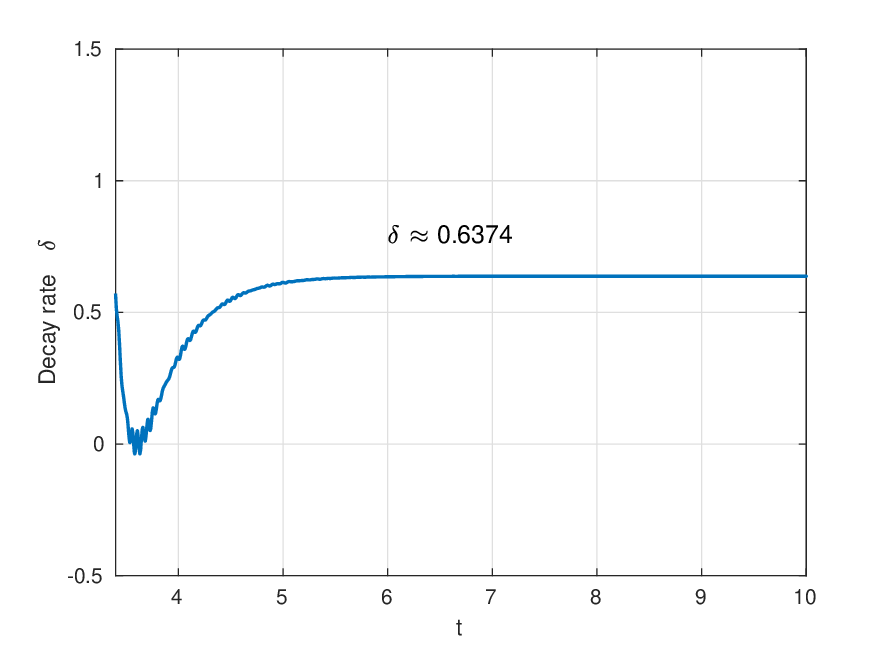}
\caption{ Estimated discrete energy decay rates for the weak-damping case in Example~\ref{example3}(ii).}
\label{fig:decaybeta0}
\end{figure}

\begin{example}\label{ex2}
We next consider $\Omega=(0,\pi)^2$ and choose
\[
\alpha=\frac{\pi^2+4}{2\pi},
\qquad
\beta=\frac{\pi}{4}.
\]
For these parameters, the homogeneous strongly damped wave
equation admits the exact solution
\[
u(x,y,t)=e^{-\pi t}\sin(x)\sin(y).
\]
This solution vanishes on $\partial\Omega$ and determines
the initial displacement and velocity as
\[
u(x,y,0)=\sin(x)\sin(y),
\qquad
u'(x,y,0)=-\pi\sin(x)\sin(y),
\quad (x,y)\in\Omega.
\]
\end{example}
\noindent
Table~\ref{errorex2} reports the errors together with the observed orders of convergence in the \(L^2\)-, \(L^\infty\)-and energy-norms. The numerical results demonstrate the accuracy of the proposed DG method and are in excellent agreement with the convergence theory developed in the preceding sections. Figures~\ref{fig:discrete_energyex2} and
\ref{fig:continuous_energyex21} show the decay of the discrete
and continuous energies, respectively.
Figure~\ref{fig:decayrate} presents the discrete energy
profiles and the estimated exponential decay rates.
\begin{table}[!ht]
\centering
\begin{tabular}{|c|cc|cc|cc|}
\hline
\multirow{2}{*}{$N$} 
 & \multicolumn{6}{c|}{$\alpha=\dfrac{\pi^{2}+4}{2\pi},\quad \beta=\dfrac{\pi}{4}$} \\ \cline{2-7}
& $\|u-U\|_{L^2}$ & Rate
& $\|u-U\|_{L^\infty}$ & Rate
& $\vertiii{u-U}$ & Rate \\
\hline
10
& 8.7198(-3) & --
& 5.6309(-3) & --
& 2.2256(-2) & -- \\
\hline
20
& 2.2219(-3) & 1.9725
& 1.4530(-3) & 1.9543
& 1.0350(-2) & 1.1046 \\
\hline
30
& 9.8762(-4) & 1.9997
& 6.4872(-4) & 1.9888
& 6.7424(-3) & 1.0570 \\
\hline
40
& 5.5493(-4) & 2.0038
& 3.6533(-4) & 1.9960
& 5.0030(-3) & 1.0372 \\
50
&3.5479(-4) & 2.0046
& 2.3389(-4) & 1.9956
& 3.9782(-3) & 1.0271 \\
\hline
\end{tabular}
\caption{Convergence history for Example~~\ref{ex2}: errors in $L^2$, $L^\infty$-norms, and energy norm, with observed rates.}
\label{errorex2}
\end{table}
\begin{example}\label{example3}
This example compares the DG approximation for purely strong damping and purely weak damping on $\Omega=(0,\pi)^2$.
The respective parameter choices are
\[
\text{(i)}\quad
\alpha=0,\qquad \beta=\frac{\pi^2+2}{2\pi},
\qquad
\text{(ii)}\quad
\alpha=\frac{\pi^2+2}{\pi},\qquad \beta=0.
\]
Both cases admit the exact solution
\[
u(x,y,t)=e^{-\pi t}\sin(x)\sin(y).
\]
We impose homogeneous Dirichlet boundary conditions and
initialise the problem with
\[
u(x,y,0)=\sin(x)\sin(y),
\qquad
u'(x,y,0)=-\pi\sin(x)\sin(y),
\quad (x,y)\in\Omega.
\]
\end{example}\noindent
Tables~\ref{errorex3alpha0} and~\ref{errorex3beta0} present the errors and the corresponding convergence rates in the $L^2$, $L^\infty$, and DG energy norms for Cases~(i) and~(ii), respectively. The observed convergence rates agree with the theoretical error bounds
for the DG approximation. Figures~\ref{fig:energy_decay_comparison} and
\ref{fig:energy_decay_comparisonex3} compare the discrete
and continuous energy histories for the purely strong
and purely weak damping cases, respectively. 
The corresponding estimated decay rates are presented in
Figures~\ref{fig:decayratealpha0} and~\ref{fig:decaybeta0}.
\begin{table}[!ht]
\centering
\begin{tabular}{|c|cc|cc|cc|}
\hline
\multirow{2}{*}{$N$}
& \multicolumn{6}{c|}{$\alpha=0,\quad \beta=\dfrac{\pi^{2}+2}{2\pi}$} \\ \cline{2-7}
& $\|u-U\|_{L^2}$ & Rate
& $\|u-U\|_{L^\infty}$ & Rate
& $\vertiii{u-U}$ & Rate \\
\hline
10
& 4.0163(-3) & --
& 6.4182(-3) & --
& 5.0388(-2) & -- \\
\hline
20
& 1.0477(-3) & 1.9386
& 1.6887(-3) & 1.9262
& 2.3754(-2) & 1.0849 \\
\hline
30
& 4.7099(-4) & 1.9719
& 7.5656(-4) & 1.9803
& 1.5514(-2) & 1.0506 \\
\hline
40
& 2.6612(-4) & 1.9844
& 4.2662(-4) & 1.9914
& 1.1522(-2) & 1.0341 \\
\hline
50
& 1.7069(-4) & 1.9901
& 2.7332(-4) & 1.9953
& 9.1658(-3) & 1.0252 \\
\hline
\end{tabular}
\caption{Errors in the $L^2$, $L^\infty$, and energy norms together with the corresponding convergence rates for Example \ref{example3}(i).}
\label{errorex3alpha0}
\end{table}
\begin{table}[!ht]
\centering
\begin{tabular}{|c|cc|cc|cc|}
\hline
\multirow{2}{*}{$N$}
& \multicolumn{6}{c|}{$\alpha=\dfrac{\pi^{2}+2}{\pi},\quad \beta=0$} \\ \cline{2-7}
& $\|u-U\|_{L^2}$ & Rate
& $\|u-U\|_{L^\infty}$ & Rate
& $\vertiii{u-U}$ & Rate \\
\hline
10
& 1.6128(-2) & --
& 1.0675(-2) & --
& 2.4924(-2) & -- \\
\hline
20
& 4.2354(-3) & 1.9290
& 2.7475(-3) & 1.9581
& 7.2256(-3) & 1.7864 \\
\hline
30
& 1.9015(-3) & 1.9751
& 1.2273(-3) & 1.9876
& 3.7668(-3) & 1.6066 \\
\hline
40
& 1.0737(-3) & 1.9866
& 6.9282(-4) & 1.9875
& 2.4499(-3) & 1.4953 \\
\hline
50
& 6.8853(-4) & 1.9913
& 4.4389(-4) & 1.9951
& 1.8062(-3) & 1.3661 \\
\hline
\end{tabular}
\caption{Errors in the $L^2$, $L^\infty$, and DG energy norms together with the corresponding convergence rates for the Example \ref{example3}(ii).}
\label{errorex3beta0}
\end{table}
    
\begin{figure}[!ht]
    \centering
    \begin{subfigure}[b]{0.48\textwidth}
        \centering
        \includegraphics[width=\textwidth, height=5.5cm, keepaspectratio=false]{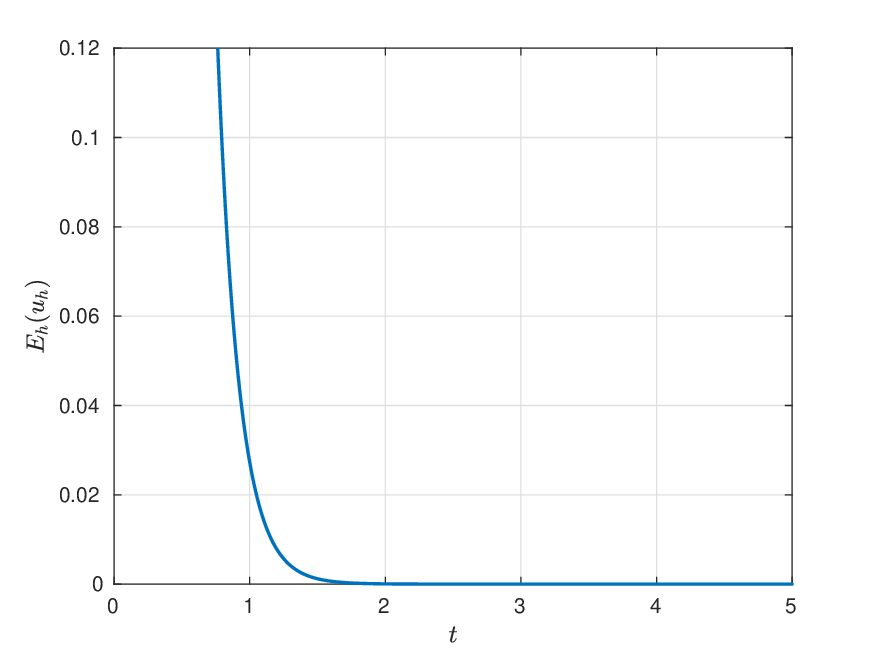}
        \caption{Discrete energy $E_h(u_h)$}
        \label{deltabenchmark1}
    \end{subfigure}
    \hfill
    \begin{subfigure}[b]{0.48\textwidth}
        \centering
        \includegraphics[width=\textwidth, height=5.5cm, keepaspectratio=false]{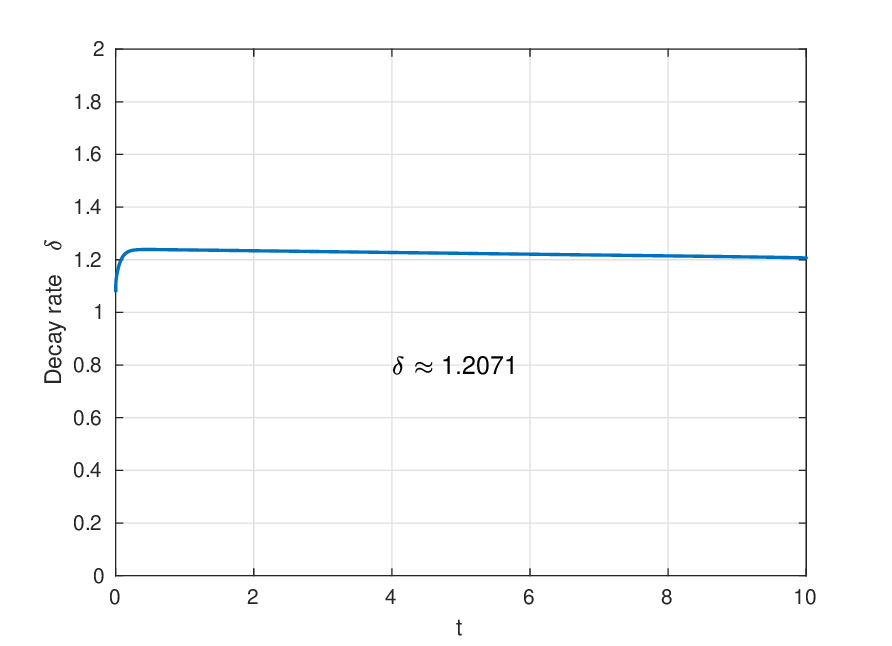}
        \caption{Energy decay with rate bound}
        \label{benchmarkdelta}
    \end{subfigure}
    \caption{Time evolution of the discrete energy $E_h(u_h)$ and comparison of estimated decay rates for Example~\ref{benchmark_ex}.}
    \label{fig:energy_benchmark}
\end{figure}

\begin{example}\label{benchmark_ex}
We examine spatial convergence on $\Omega=(0,\pi)^2$
for the strongly damped wave equation with zero forcing,
using
\[
\alpha=\frac{\pi^2+4}{2\pi},
\qquad \beta=\frac{\pi}{4}.
\]
Homogeneous Dirichlet conditions are imposed on
$\partial\Omega\times(0,T]$, and the initial data are
\[
u(x,y,0)=\sin(x)\sin(y),
\qquad
u'(x,y,0)=-\sin(x)\sin(y).
\]
A reference solution is computed using discontinuous
$P_1$ elements on a fine mesh with $N=80$.
Errors on the coarser meshes are measured relative to
this reference solution and used to calculate the
observed convergence rates.
\end{example}
Table~\ref{benchmark_tab} shows convergence rates approaching
two in the $L^2$ and $L^\infty$ norms and one in the energy norm.
Figure~\ref{fig:energy_benchmark} presents the discrete energy
history and its comparison with the exponential decay bound.
\begin{table}[!ht]
\centering
\begin{tabular}{|c|cc|cc|cc|}
\hline
\multirow{2}{*}{$N$} 
 & \multicolumn{6}{c|}{$\alpha=\dfrac{\pi^{2}+4}{2\pi},\quad \beta=\dfrac{\pi}{4}$} \\ \cline{2-7}
& $\|u-U\|_{L^2}$ & Rate
& $\|u-U\|_{L^\infty}$ & Rate
& $\vertiii{u-U}$ & Rate \\
\hline
10
& 1.0015(-1) & --
& 7.7097(-2) & --
& 9.8116(-1) & -- \\
\hline
20
& 2.9650(-2) & 1.7560
& 2.3231(-2) & 1.7306
& 5.1162(-1) & 0.9394 \\
\hline
30
& 1.4024(-2) & 1.8465
& 1.0846(-2) & 1.8786
& 3.4426(-1) & 0.9771 \\
\hline
40
& 8.0369(-3) & 1.9351
& 6.3249(-3) & 1.8747
& 2.5858(-1) & 0.9948 \\
\hline
50
& 5.0678(-3) & 2.0665
& 3.8945(-3) & 2.1731
& 2.0666(-1) & 1.0044 \\
\hline
\end{tabular}
\caption{Errors relative to the reference solution in the $L^2$, $L^\infty$, and energy norms, with observed convergence rates, for Example~\ref{benchmark_ex}.}
\label{benchmark_tab}
\end{table}
\subsection{DG Method for Fourth Order Evolution Equations}
In the remaining part, we discuss  DG scheme for the fourth order evolution equations \eqref{bihar1}-\eqref{bihar2}  given in Section 5. All computations in this part are carried out using
\texttt{Firedrake}.

\begin{example}\label{ex1bihar}
Consider $\Omega=(0,1)^2$ and choose
$\alpha=\pi$ and $\beta=1/\pi$.
The exact solution is
\[
u(x,y,t)=e^{-\pi t}\sin^2(\pi x)\sin^2(\pi y),
\]
with initial displacement and velocity
\[
u(x,y,0)=\sin^2(\pi x)\sin^2(\pi y),
\qquad
u'(x,y,0)=-\pi\sin^2(\pi x)\sin^2(\pi y).
\]
Both the solution and its normal derivative vanish on
the boundary, so the homogeneous clamped conditions
\[
u=0,\qquad \frac{\partial u}{\partial n}=0
\quad\text{on }\partial\Omega\times(0,\infty)
\]
are satisfied.
\end{example}
Tables~\ref{tab:convergence_rates1} and
\ref{tab:convergence_ratesP3} report the errors and observed
convergence rates for the $P_2$ and $P_3$ DG  approximations,
respectively. The $L^2$ and $L^\infty$ errors exhibit approximately second-order convergence for $P_2$ and fourth-order convergence
for $P_3$, while the DG energy errors converge at first and
second order, respectively, matching the predicted orders
of convergence. Tables~\ref{tab:convergence_rates2} and
\ref{tab:convergence_rates1P3COIP} present the convergence
results for the $P_2$ and $P_3$ $C^0$ interior penalty
approximations, respectively. As the mesh is refined,
the $L^2$ and $L^\infty$ convergence rates approach two
for $P_2$ and four for $P_3$, while the energy-norm rates
approach one and two, respectively, in agreement with
the predicted orders of convergence.
\begin{table}[h!]
\centering
\begin{tabular}{|c|cc|cc|cc|}
\hline
\multirow{2}{*}{$N$} 
& \multicolumn{6}{c|}{$\alpha=\pi,\quad \beta=\frac{1}{\pi}$} \\ \cline{2-7}
& $\|u-U\|_{L^2}$ & Rate 
& $\|u-U\|_{L^\infty}$ & Rate 
& $\vertiii{u-U}_{dG}$ & Rate \\
\hline
10  
& 4.0610(-3) & --    
& 1.0820(-2) & --    
& 3.5430(-1) & -- \\
\hline
20 
& 1.0140(-3) & 2.0019 
& 2.6990(-3) & 2.0036 
& 1.8050(-1) & 0.9731 \\
\hline
30 
& 4.3790(-4) & 2.0705 
& 1.1630(-3) & 2.0760 
& 1.2130(-1) & 0.9804 \\
\hline
40 
& 2.3970(-4) & 2.0945 
& 6.3600(-4) & 2.0985 
& 9.1290(-2) & 0.9871 \\
\hline
50 
& 1.5180(-4) & 2.0484 
& 4.0220(-4) & 2.0535 
& 7.3190(-2) & 0.9905 \\
\hline
\end{tabular}
\caption{Convergence history for $P_2$ spatial discretization: $L^2$, $L^\infty$, and dG errors for dG FEM method with $\sigma_e=10,\,\sigma_{dG}^1=10,\,\sigma_{dG}^2=15$.}
\label{tab:convergence_rates1}
\end{table}
\begin{table}[h!]
\centering
\begin{tabular}{|c|cc|cc|cc|}
\hline
\multirow{2}{*}{$N$} 
& \multicolumn{6}{c|}{$\alpha=\pi,\quad \beta=\frac{1}{\pi}$} \\ \cline{2-7}
& $\|u-U\|_{L^2}$ & Rate 
& $\|u-U\|_{L^\infty}$ & Rate 
& $\vertiii{u-U}_{dG}$ & Rate \\
\hline
10  
& 7.1540(-5) & --    
& 1.9560(-4) & --    
& 2.8690(-2) & -- \\
\hline
20 
& 4.4530(-6) & 4.0060 
& 1.2220(-5) & 4.0004 
& 6.6350(-3) & 2.1124 \\
\hline
30 
& 8.7820(-7) & 4.0038 
& 2.4150(-6) & 3.9989 
& 2.8650(-3) & 2.0707 \\
\hline
40 
& 2.7950(-7) & 3.9798 
& 7.6890(-7) & 3.9785 
& 1.5880(-3) & 2.0514 \\
\hline
50 
& 1.1480(-7) & 3.9863 
& 3.1600(-7) & 3.9855 
& 1.0070(-3) & 2.0398 \\
\hline
\end{tabular}
\caption{Convergence history for $P_3$ spatial discretization: $L^2$, $L^\infty$, and dG errors for dG FEM method with $\sigma_e=200,\,\sigma_{dG}^1=300,\,\sigma_{dG}^2=300$.}
\label{tab:convergence_ratesP3}
\end{table}

\begin{table}[h!]
\centering
\begin{tabular}{|c|cc|cc|cc|}
\hline
\multirow{2}{*}{$N$} 
& \multicolumn{6}{c|}{$\alpha=\pi,\quad \beta=\frac{1}{\pi}$} \\ \cline{2-7}
& $\|u-U\|_{L^2}$ & Rate 
& $\|u-U\|_{L^\infty}$ & Rate 
& $\vertiii{u-U}_{dG}$ & Rate \\
\hline
10  
& 3.3540(-3) & --    
& 5.4320(-3) & --    
& 2.1500(-1) & -- \\
\hline
20 
& 8.8770(-4) & 1.9176 
& 1.5250(-3) & 1.8323 
& 1.0450(-1) & 1.0415 \\
\hline
30 
& 3.9970(-4) & 1.9681 
& 7.0300(-4) & 1.9108 
& 6.8580(-2) & 1.0378 \\
\hline
40 
& 2.2590(-4) & 1.9829 
& 4.0160(-4) & 1.9458 
& 5.1060(-2) & 1.0252 \\
\hline
50 
& 1.4490(-4) & 1.9892 
& 2.5920(-4) & 1.9625 
& 4.0690(-2) & 1.0177 \\
\hline
\end{tabular}\caption{Convergence history for $P_2$ spatial discretization: $L^2$, $L^\infty$, and dG errors for $C^0$-IP method with $\sigma_{IP}=10$.}
\label{tab:convergence_rates2}
\end{table}
\begin{table}[h!]
\centering
\begin{tabular}{|c|cc|cc|cc|}
\hline
\multirow{2}{*}{$N$} 
& \multicolumn{6}{c|}{$\alpha=\pi,\quad \beta=\frac{1}{\pi}$} \\ \cline{2-7}
& $\|u-U\|_{L^2}$ & Rate 
& $\|u-U\|_{L^\infty}$ & Rate 
& $\vertiii{u-U}_{dG}$ & Rate \\
\hline
10  
& 1.9710(-2) & --    
& 1.9910(-2) & --    
& 2.5230(-2) & -- \\
\hline
20 
& 1.9650(-3) & 3.3264 
& 1.9780(-3) & 3.3313 
& 6.2550(-3) & 2.0117 \\
\hline
30 
& 3.9610(-4) & 3.9501 
& 3.9870(-4) & 3.9498 
& 2.7860(-3) & 1.9950 \\
\hline
40 
& 1.2570(-4) & 3.9883 
& 1.2660(-4) & 3.9878 
& 1.5680(-3) & 1.9973 \\
\hline
50 
& 5.1550(-5) & 3.9959 
& 5.1910(-5) & 3.9954 
& 1.0040(-3) & 1.9987 \\
\hline
\end{tabular}
\caption{Convergence history for $P_3$ spatial discretization: $L^2$, $L^\infty$, and dG errors for $C^0$-IP method with $\sigma_{IP}=20$.}
\label{tab:convergence_rates1P3COIP}
\end{table}

\begin{example}\label{ex2bihar_bench}
We next examine convergence for the homogeneous fourth-order damped wave problem on $\Omega=(0,1)^2$, with
$\alpha=\pi$ and $\beta=1/\pi$.
The boundary conditions are
\[
u=0,\qquad \frac{\partial u}{\partial n}=0
\quad\text{on }\partial\Omega\times(0,T].
\]
The prescribed initial displacement and velocity are
\[
u(x,y,0)=\sin^2(\pi x)\sin^2(\pi y),
\qquad
u'(x,y,0)=-\sin^2(\pi x)\sin^2(\pi y).
\]
For the DG computations, a reference solution is obtained
using discontinuous $P_2$ elements on a uniform mesh with
$N=120$. Errors on the test meshes are measured relative
to this reference solution, and convergence rates are
calculated from successive refinements.
\end{example}

\begin{figure}[!htbp]
\centering

\includegraphics[
  width=0.65\linewidth,
  trim={0bp 0bp 350bp 0bp},
  clip
]{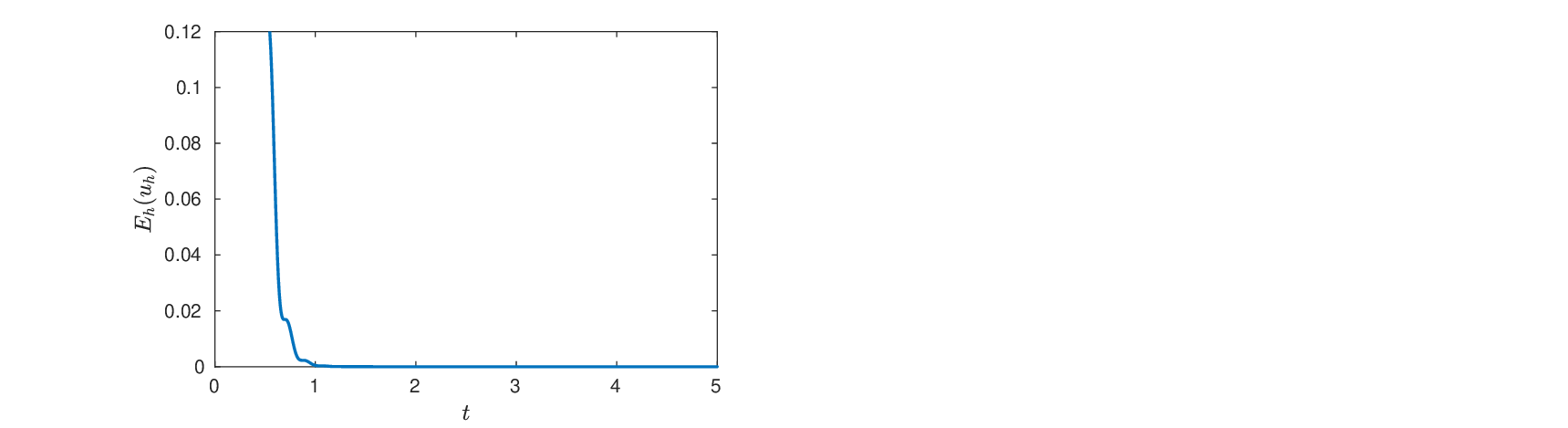}
\caption{Discrete energy $E_h(u_h)$ versus time
for Example~\ref{ex2bihar_bench}.}
\label{energy_benchmark_dGFEM}
\end{figure}

Tables~\ref{tab:convergence_rates_ref} and
\ref{benchmark_dGFEM} report the results for the $P_2$ DG
and $C^0$ interior penalty methods, respectively.
The $L^2$ and $L^\infty$ errors exhibit approximately
second-order convergence, although the DG rates vary
across refinement levels. The energy-norm rates are
close to one on the finer DG meshes and decrease towards
one for the $C^0$ interior penalty method.
Figure~\ref{energy_benchmark_dGFEM} shows the evolution
of the discrete energy.\begin{table}[h!]
\centering
\begin{tabular}{|c|cc|cc|cc|}
\hline
\multirow{2}{*}{$N$} 
& \multicolumn{6}{c|}{$\alpha=\pi,\quad \beta=\frac{1}{\pi}$} \\ \cline{2-7}
& $\|u-U\|_{L^2}$ & Rate 
& $\|u-U\|_{L^\infty}$ & Rate 
& $\vertiii{u-U}_{dG}$ & Rate \\
\hline
10  
& 1.0028(-4) & --    
& 2.7895(-4) & --    
& 1.6279(-2) & -- \\
\hline
20 
& 2.3301(-5) & 2.1056 
& 4.8106(-5) & 2.5357 
& 1.0917(-2) & 0.5764 \\
\hline
30 
& 1.1786(-5) & 1.6811 
& 2.5474(-5) & 1.5680 
& 6.7977(-3) & 1.1684 \\
\hline
40 
& 6.3021(-6) & 2.1760 
& 1.3724(-5) & 2.1501 
& 4.8564(-3) & 1.1690 \\
\hline
\end{tabular}
\caption{Errors relative to the fine-mesh reference solution
and observed convergence rates for the $P_2$ DG method
in Example~\ref{ex2bihar_bench}.}
\label{tab:convergence_rates_ref}
\end{table}
\begin{table}[!ht]
\centering
\begin{tabular}{|c|cc|cc|cc|}
\hline
\multirow{2}{*}{$N$} 
& \multicolumn{6}{c|}{$\alpha=\pi,\quad \beta=\frac{1}{\pi}$} \\ \cline{2-7}
& $\|u-U\|_{L^2}$ & Rate 
& $\|u-U\|_{L^\infty}$ & Rate 
& $\vertiii{u-U}_{dG}$ & Rate \\
\hline
10  
& 6.4904(-4) & --    
& 1.5886(-3) & --    
& 2.3201(-2) & -- \\
\hline
20 
& 1.6094(-4) & 2.0118 
& 3.8892(-4) & 2.0302 
& 6.6015(-3) & 1.8134 \\
\hline
30 
& 6.9629(-5) & 2.0664 
& 1.6753(-4) & 2.0771 
& 3.4311(-3) & 1.6140 \\
\hline
40 
& 3.8353(-5) & 2.0730 
& 9.2108(-5) & 2.0794 
& 2.2597(-3) & 1.4518 \\
\hline
50 
& 2.3980(-5) & 2.1044 
& 5.7538(-5) & 2.1086 
& 1.6950(-3) & 1.2885 \\
\hline
\end{tabular}
\caption{Time evolution of the discrete energy $E_h(u_h)$
for Example~\ref{ex2bihar_bench}.}
\label{benchmark_dGFEM}
\end{table}

\section{Summary}\label{s7}
This paper develops and analyses discontinuous Galerkin
approximations for wave equations with weak and strong damping. Energy arguments establish exponential decay for both the semidiscrete and fully discrete solutions, with decay rates independent of the discretisation parameters under the stated
assumptions and time-step restrictions. The analysis also provides error estimates in the energy and $L^2$ norms that account for the long-time decay of the exact solution. Extensions to
nonhomogeneous problems and fourth-order damped wave equations are considered, with the latter treated using symmetric interior penalty DG and $C^0$ interior penalty methods. Numerical experiments
support the theoretical convergence estimates, illustrate the
decay of the discrete energy, and examine the influence of the
damping parameters on the computed solutions.
\section*{Statements and Declarations}
\subsection*{Ethical Approval} Not Applicable
\subsection*{Data Availability}
The codes during the current study are available from the corresponding author on reasonable request.
\subsection*{Conflict of Interest}
The authors declare that they have no conflict of interest.
\subsection*{Funding}
The authors acknowledge financial support from ANRF, the
Government of India, via Project No.~CRG/2023/000721.

\subsection*{Author Contributions}
 Amiya K. Pani, P. Danumjaya and Anil Kumar jointly conceptualized the project, supervised, developed the analysis, wrote, and revised the manuscript. Devendra K. Swain also participated in conceptualizing the project, developed the analysis, validated the results, and contributed to the writing the manuscript. All authors contributed equally to the preparation of this manuscript and  approved the final manuscript.

\subsection*{Acknowledgments} 
The first author acknowledges the financial support from ANRF project No. ~CRG/2023/000721.

\bibliographystyle{plain}
\bibliography{DG-references-uniform}
\end{document}